\documentclass[reqno,a4paper,12pt]{amsart}
\usepackage{amssymb}
\usepackage[dvips, lmargin=3cm, rmargin=3cm, tmargin=4cm, bmargin=4cm]{geometry}
\usepackage[colorlinks = true,
            linkcolor = blue,
            urlcolor  = blue,
            citecolor = red,
            anchorcolor = blue]{hyperref}
\usepackage{hyperref}
\usepackage[T1]{fontenc}
\usepackage{enumitem}
\usepackage{dsfont}
\everymath{\displaystyle}
\newtheorem{theorem}{Theorem}[section]
\newtheorem{definition}{Definition}[section]
\newtheorem{remark}{Remark}[section]
\newtheorem{lemma}{Lemma}[section]
\newtheorem{proposition}{Proposition}[section]
\newtheorem{corollary}{Corollary}[section]
\newtheorem{example}{Example}

\numberwithin{equation}{section}
\begin{document}
\title[Generalized  fractional Sobolev Space  ]{General fractional Sobolev Space with variable exponent and applications to nonlocal problems}
\author[E. Azroul, A. Benkirane and  M. Shimi]
{E. Azroul$^1$, A. Benkirane$^2$ and  M. Shimi$^3$}
\address{E. Azroul, A. Benkirane and  M. Shimi \newline
 Sidi Mohamed Ben Abdellah
 University,
 Faculty of Sciences Dhar Al Mahraz, Laboratory of Mathematical Analysis and Applications, Fez, Morocco.}
\email{$^1$elhoussine.azroul@gmail.com}
\email{$^2$abd.benkirane@gmail.com}
\email{$^3$mohammed.shimi2@usmba.ac.ma}
\maketitle

\vspace{0.2cm}
\begin{abstract}
In this paper, we extend the fractional Sobolev spaces with variable exponents $W^{s,p(x,y)}$ to include the general fractional case $W^{K,p(x,y)}$, where $p$ is a variable exponent, $s\in (0,1)$ and $K$ is a suitable kernel. We are concerned with some qualitative properties of the space $W^{K,p(x,y)}$ (completeness, reflexivity, separability and density). Moreover,  we prove a continuous and compact embedding theorem of these spaces into variable exponent Lebesgue spaces. As  applications,
 we discus the existence of a nontrivial solution for a nonlocal $p(x,.)$-Kirchhoff type problem. Further, we establish the existence and uniqueness of a solution for a variational  problem involving the integro-differential operator of elliptic type $\mathcal{L}^{p(x,.)}_K$. 
\end{abstract}

\section{Introduction}

Our main goal in this paper is to extend the fractional Sobolev spaces with variable exponents to cover the nonlocal general  case with singular kernel. For this,
we begin this work by remembering the definition of  fractional Sobolev  spaces with variable exponent, see for instance  \cite{SM, SH, 8, 12}.

Let $\Omega$ be a smooth bounded open set in $\mathbb{R}^{N}$. We start by fixing $s\in (0,1)$ and let $p:\overline{\Omega}\times\overline{\Omega}\longrightarrow(1,+\infty)$ be a continuous bounded function. We assume that
\begin{equation}\label{1}
1<p^{-}=\underset{(x,y)\in \overline{\Omega}\times\overline{\Omega}}{\min}p(x,y)\leqslant p(x,y)\leqslant p^{+}=\underset{(x,y)\in \overline{\Omega}\times\overline{\Omega}}{\max}p(x,y)<+\infty
\end{equation} 
and
\begin{equation}\label{2}
p ~~~\text{is symmetric, that is,}~~ p(x,y)=p(y,x)\quad\text{for all}\quad(x,y)\in\overline{\Omega}\times\overline{\Omega}.
\end{equation} 
Let denote by :
$$\bar{p}(x)=p(x,x) \quad\text{ for all }\quad x\in \overline{\Omega}.$$

  We define the fractional Sobolev space with variable exponent via the Gagliardo approach as follows,\\ $$\hspace{-11.2cm}W=W^{s,p(x,y)}(\Omega)$$
  $$\hspace{0.7cm}=\bigg\{u\in L^{\bar{p}(x)}(\Omega): \int_{\Omega\times\Omega}\frac{|u(x)-u(y)|^{p(x,y)}}{\lambda^{p(x,y)}|x-y|^{sp(x,y)+N}}~dxdy <+\infty,~~ \text{for some}~~\lambda>0\bigg\},$$
 where $L^{\bar{p}(x)}(\Omega)$ is the Lebesgue space with variable exponent, (see Section $\ref{200}$).\\
The space $W^{s,p(x,y)}(\Omega)$ is a Banach space (see \cite{8}) if it is endowed with the norm, $$\|u\|_{W^{s,p(x,y)}(\Omega)}=\|u\|_{W}=\|u\|_{L^{\bar{p}(x)}(\Omega)}+[u]_{s,p(x,y)},$$
where $[.]_{s,p(x,y)}$ is a Gagliardo semi-norm with variable exponent, which is defined
by $$[u]_{s,p(x,y)}= \inf \bigg\{\lambda>0:\int_{\Omega\times\Omega}\frac{|u(x)-u(y)|^{p(x,y)}}{\lambda^{p(x,y)}|x-y|^{sp(x,y)+N}}~dxdy \leqslant1 \bigg\}.$$ 
The space $(W,\|.\|_{W})$ is separable and  reflexive, see (\cite[Lemma 3.1]{2}).\\

Let us consider the fractional $p(x,.)$-Laplacian operator given by 
$$(-\Delta_{p(x,.)})^{s}u(x)=p.v.\int_{\Omega}\frac{|u(x)-u(y)|^{p(x,y)-2}(u(x)-u(y))}{|x-y|^{n+sp(x,y)}}~dy~,~~~~~ \text{for all } x \in \Omega$$
where $p.v.$ is a commonly used abbreviation in the principal value sense.
One typical feature of this operator  is the nonlocality, in the sense that the value of $(-\Delta_{p(x,.)})^{s}u(x)$ at any point $x\in \Omega$ depends not only on the values of $u$ on $\Omega$, but actually on the entire space $\mathbb{R}^{N}$.\\

Note that the operator $(-\Delta_{p(x,.)})^{s}$ is the fractional version of well known $p(x)$-Laplacian operator $\Delta_{p(x)}u(x)= div\big(|\nabla u(x)|^{p(x)-2}u(x)\big)$. On the other hand, we remark that in the constant exponent case it is know as the fractional $p$-Laplacian operator $(-\Delta)_{p}^{s}$.  This nonlinear operator is consistent, up to some normalization constant depending upon $N$ and $s$, with the linear fractional Laplacian $(-\Delta)^{s}$ in the case $p=2$. The interest for this last operator and more generally pseudo-differential operators, has constantly increased over the last few years, although such operators have been a classical topic of functional analysis since long ago. Nonlocal operators such as $(-\Delta)^{s}$ and its generalisation $\mathcal{L}_K$ (see for instance \cite{11, 13, 17, 18, 19}) naturally arise in continuum mechanics, phase transition phenomena, population dynamics and game theory, as they are the typical outcome of stochastical stabilization of L\'{e}vy processes, see e.g. \cite{6, 15, 16}. We refer the reader to \cite{4,9} and to the references included for a self-contained overview of the basic properties of fractional Sobolev spaces and fractional Laplacian operator.\\

Now, we introduce the nonlocal integro-differential operator of elliptic type  $\mathcal{L}^{p(x,.)}_K$ that generalizes the operator $(-\Delta_{p(x,.)})^{s}$,  as follows

{\small$$
       \mathcal{L}^{p(x,.)}_K(u(x))=p.v.\int_{\mathbb{R}^N}|u(x)-u(y)|^{p(x,y)-2}(u(x)-u(y))K(x,y)~dy,~~~~~ \text{for all } x \in \mathbb{R}^N$$
       $$
    \hspace{2.6cm}  =\lim\limits_{\varepsilon\longrightarrow0}\int_{\mathbb{R}^N\setminus B_{\varepsilon}(x)}|u(x)-u(y)|^{p(x,y)-2}(u(x)-u(y))K(x,y)~dy,~~~~~ \text{for all } x \in \mathbb{R}^N
     .$$} where $p.v.$ is a commonly used abbreviation in the principal value sense,  $p:\mathbb{R}^N \times\mathbb{R}^N\longrightarrow(1,+\infty)$ is a continuous bounded function satisfy (\ref{1}), (\ref{2}) and \begin{equation}\label{3}
      p((x,y)-(z,z))=p(x,y),\hspace{0.5cm}~~~~~ \text{for all }(x,y),(z,z)\in  \mathbb{R}^N \times\mathbb{R}^N .
      \end{equation}
  The kernel $K:\mathbb{R}^N \times\mathbb{R}^N\longrightarrow (0,+\infty)$ is a measurable function with the following properties:
   \begin{equation}\label{4}
 K(x,y)=K(y,x)\hspace{0.3cm}\text{for any } (x,y)\in  \mathbb{R}^N \times\mathbb{R}^N ,
  \end{equation}  
there exists $ k_0 >0$  such that 
 \begin{equation}\label{5}
  K(x,y)\geqslant k_0|x-y|^{-(N+sp(x,y))},~ \text{for any } (x,y)\in  \mathbb{R}^N \times\mathbb{R}^N~ and~x\neq y,
  \end{equation} 
 
  \begin{equation}\label{6}
  m K\in L^{1}( \mathbb{R}^N \times\mathbb{R}^N), \textit{where}~ m(x,y)=\min\big\{ 1, |x-y|^{p(x,y)}\big\}. 
  \end{equation}
  
  A typical example for $K$  is given by singular kernel $K(x,y)=|x-y|^{-(N+sp(x,y))}$. In this case  $ \mathcal{L}^{p(x,.)}_K=(-\Delta_{p(x,.)})^{s}$.
  An other example for $K$ is given by the kernel $$K_1(x,y)= |x-y|^{-(N+sp(x,y))}a(x-y),$$
  where  $a:\mathbb{R}^{N} \longrightarrow[1,+\infty)$ is a bounded function, that is, $a\in L^{\infty}(\mathbb{R}^{N})$. It is easy to see that $K_1$ satisfy the assumptions (\ref{4})-(\ref{6}).\\
  
   This paper is organized as follows. In Section $\ref{200}$, we give some definitions and fundamental properties of the spaces $L^{q(x)}$ and $W^{s,p(x,y)}$.  In Section $\ref{300}$, we compare the space $W^{s,p(x,y)}$ with $W^{K,p(x,y)}$ and  we study the completeness, reflexivity, separability, and density of these spaces. Moreover, we prove a continuous and compact embedding theorem of these spaces into variable exponent Lebesgue spaces. In section \ref{400}, we prove some basic properties of the operator $\mathcal{L}^{p(x.,)}_K$. As  applications, in Section \ref{500}, we show the existence of a nontrivial solution for a nonlocal $p(x,.)$-Kirchhoff type problem by means of mountain pass theorem. Finally,  we apply the Minty-Browder theorem
   to establish the existence and uniqueness of a solution for a variational problem involving the the integro-differential operator $\mathcal{L}^{p(x,.)}_K$. 
  
  \section{Some preliminary results}\label{200}
   In this section, we recall some necessary properties of variable exponent spaces. For more details we refer the reader to \cite{7,10,14}, and the references therein.\\ 
    Consider the set
     $$C_+(\overline{\Omega})=\big\{q\in C(\overline{\Omega}): q(x)>1,~~ \text{for all} x \in\overline{\Omega}\big\}.$$
     For all $q\in C_+(\overline{\Omega}) $, we define $$q^{+}= \underset{x\in \overline{\Omega}}{\sup}~q(x) ~~~\text{and}~~~ q^{-}= \underset{x\in \overline{\Omega}}{\inf}~q(x).$$
  Such that
   \begin{equation}\label{7}
          1<q^{-}\leqslant q(x)\leqslant q^{+}<+\infty.
          \end{equation} 
  For any  $q\in C_+(\overline{\Omega}) $, we define the variable exponent Lebesgue space as $$L^{q(x)}(\Omega)=\bigg\{u:\Omega\longrightarrow \mathbb{R} ~~\text{measurable}: \int_{\Omega}|u(x)|^{q(x)}dx<+\infty
  \bigg\}.$$
  
  This vector space endowed with the \textit{Luxemburg norm}, which is defined by
  $$\|u\|_{L^{q(x)}(\Omega)}= \inf \bigg\{\lambda>0:\int_{\Omega}\bigg|\frac{u(x)}{\lambda}\bigg|^{q(x)}dx \leqslant1 \bigg\}$$
  is a separable reflexive Banach space.\\
  Let $ \hat{q} \in C_+(\overline{\Omega})$  be the conjugate exponent of $q$, that is, $\frac{1}{q(x)}+\frac{1}{\hat{q}(x)}=1$. Then we have the following H\"{o}lder-type inequality 
  \begin{lemma}$\label{2.1}$\textnormal{(H\"{o}lder inequality)}.
  If $u\in L^{q(x)}(\Omega)$ and $v\in L^{\hat{q}(x)}(\Omega)$, so
  $$\bigg|\int_{\Omega}u v dx\bigg|\leqslant \bigg(\frac{1}{q^{-}}+\frac{1}{\hat{q}^{-}}\bigg)\|u\|_{L^{q(x)}(\Omega)}\|v\|_{L^{\hat{q}(x)}(\Omega)}\leqslant 2\|u\|_{L^{q(x)}(\Omega)}\|v\|_{L^{\hat{q}(x)}(\Omega)}.$$
  \end{lemma}
   A very important role in manipulating the generalized Lebesgue spaces with variable exponent is played by the modular of the $L^{q(x)}(\Omega)$ space, which defined by
   $$\begin{array}{clc}
   \hspace{-0.8cm}\rho_{q(.)}: L^{q(x)}(\Omega)\longrightarrow\mathbb{R}\\
    \hspace{5.6cm}u\longrightarrow\rho_{q(.)}(u)=\int_{\Omega}|u(x)|^{q(x)}dx
   \end{array}$$
  \begin{proposition}$\label{2.2.1}$
  Let $u\in  L^{q(x)}(\Omega) $, then we have,
  \begin{enumerate}[label=(\roman*)]
  \item $\|u\|_{L^{q(x)}(\Omega)}<1$ $($resp $=1$, $>1$$)$ $\Leftrightarrow$ $ \rho_{q(.)}(u)<1$ $($resp $=1$, $>1$$)$,
  \item  $\|u\|_{L^{q(x)}(\Omega)}<1$ $\Rightarrow$ $\|u\|^{q{+}}_{L^{q(x)}(\Omega)}\leqslant \rho_{q(.)}(u)\leqslant \|u\|^{q{-}}_{L^{q(x)}(\Omega)}$,
  \item  $\|u\|_{L^{q(x)}(\Omega)}>1$ $\Rightarrow$ $\|u\|^{q{-}}_{L^{q(x)}(\Omega)}\leqslant \rho_{q(.)}(u)\leqslant \|u\|^{q{+}}_{L^{q(x)}(\Omega)}$.
  \end{enumerate}
  \end{proposition}
  
  \begin{proposition}$\label{2.2.2}$
  If $u, u_k \in L^{q(x)}(\Omega) $ and $k\in \mathbb{N}$, then the following assertions are equivalent
  \begin{enumerate}[label=(\roman*)]
  \item $\lim\limits_{k\rightarrow+\infty} \|u_k-u\|_{L^{q(x)}(\Omega)}=0$,
  \item$\lim\limits_{k\rightarrow+\infty} \rho_{q(.)}(u_k-u)=0$,
  \item $u_k\longrightarrow u$ in measure in $\Omega$ and $\lim\limits_{k\rightarrow+\infty} \rho_{q(.)}(u_k)=\rho_{q(.)}(u)$.
  \end{enumerate}
  \end{proposition}
  From Theorems 1.6, 1.8 and 1.10 of \cite{10}, we obtain the following proposition:
  
  \begin{proposition}\label{2.3.1}
  Suppose that $(\ref{7})$ is satisfied. If $\Omega$ is a bounded open domain, $(L^{q(x)}(\Omega), \|u\|_{L^{q(x)}(\Omega)})$ is a reflexive uniformly convex and separable Banach space.
  \end{proposition}
  \begin{definition}
   Let $p:\overline{\Omega}\times\overline{\Omega}\longrightarrow(1,+\infty)$,
   be a continuous variable exponent and $s\in (0,1)$.
  For any $u\in W$, we define the modular\\ $\rho_{p(.,.)}:W\longrightarrow\mathbb{R}$, by 
   $$ \rho_{p(.,.)}(u)=\int_{\Omega\times\Omega}\frac{|u(x)-u(y)|^{p(x,y)}}{|x-y|^{N+sp(x,y)}}~dxdy+\int_{\Omega}|u(x)|^{\bar{p}(x)}dx
   $$
   and
    $$\|u\|_{\rho_{p(.,.)}}= \inf \bigg\{\lambda>0:\rho_{p(.,.)}\bigg(\frac{u}{\lambda}\bigg) \leqslant1 \bigg\}.$$ 
  \end{definition}
  \begin{remark}$\label{2.5}$\text{}
  \begin{enumerate}[label=(\roman*)]
  \item It is easy to see that $\|.\|_{\rho_{p(.,.)}}$ is a norm on $W$ which is equivalent to the norm $\|.\|_{W}$.
  \item $\rho_{p(.,.)}$  also check the results of Propositions $\ref{2.2.1}$ and $\ref{2.2.2}$.
  \end{enumerate}
  \end{remark}
  The following inequality can be easily obtained.
  \begin{lemma}\label{2.2}
  Let $p:\overline{\Omega}\times\overline{\Omega}\longrightarrow(1,+\infty)$,
   be a continuous variable exponent and $s\in (0,1)$.
  For any $u\in W_0$, we have
  \begin{enumerate}[label=(\roman*)]
  {\small \item $1\leqslant [u]_{s,p(x,y)}$ $\Rightarrow$ $[u]^{p^{-}}_{s,p(x,y)}\leqslant\int_{\Omega\times\Omega}\frac{|u(x)-u(y)|^{p(x,y)}}{|x-y|^{N+sp(x,y)}}~dxdy\leqslant [u]^{p^{+}}_{s,p(x,y)}$,}
  {\small\item  $ [u]_{s,p(x,y)}\leqslant1$ $\Rightarrow$ $[u]^{p^{+}}_{s,p(x,y)}\leqslant\int_{\Omega\times\Omega}\frac{|u(x)-u(y)|^{p(x,y)}}{|x-y|^{N+sp(x,y)}}~dxdy\leqslant [u]^{p^{-}}_{s,p(x,y)}$.}
  \end{enumerate}
  \end{lemma}
  In \cite{12}, the authors introduce the variable exponent  Sobolev fractional space as follows 
   $$\hspace{-10.6cm} E= W^{s,q(x),p(x,y)}(\Omega)$$
    $$\hspace{0.5cm}=\bigg\{u\in L^{q(x)}(\Omega): \int_{\Omega\times\Omega}\frac{|u(x)-u(y)|^{p(x,y)}}{\lambda^{p(x,y)}|x-y|^{sp(x,y)+N}}~dxdy <+\infty,~~ \text{for some}~~\lambda>0\bigg\},$$
  where $q:\overline{\Omega}\longrightarrow(1,+\infty)$ be a continuous function satisfying  (\ref{7}).
  
   We would like to mention that the continuous and compact embedding theorem has been proved in \cite{12} under the assumption {\small$q(x) > \bar{p}(x)=p(x,x)$}. The authors in \cite{SH} give a slightly different version of continuous compact embedding theorem assuming that $q(x)=\bar{p}(x)=p(x,x)$..
  \begin{theorem}$\label{1.1}$\textnormal{(\cite{SH})}.
  Let $\Omega$ be a Lipschitz bounded  domain in $\mathbb{R}^{N}$ and let $s\in (0,1)$. Let $p:\overline{\Omega}\times\overline{\Omega}\longrightarrow(1,+\infty)$
   be a continuous function satisfies $(\ref{1})$ and  $(\ref{2})$ with $s p^{+}<N$ 
   . Let $r:\overline{\Omega}\longrightarrow(1,+\infty)$  be a continuous variable exponent such that  $$1< r^{-}=\underset{x\in \overline{\Omega}}{\min}~r(x)\leqslant r(x) <p^{\ast}_{s}(x)= \frac{N\bar{p}(x)}{N-s\bar{p}(x)} ~~ \text{ for all } x\in\overline{\Omega}.$$
    Then, there exists a constant $C=C(N,s,p,r, \Omega)>0$ such that, for any $u\in W$,$$\|u\|_{L^{r(x)}(\Omega)}\leqslant C\|u\|_{W}.$$
   That is, the space $W$ is continuously embedded in $L^{r(x)}(\Omega)$. Moreover, this
    embedding is compact.
  \end{theorem}
  
  \begin{remark}$\label{1.2}$Let $W_{0}$ denotes the closure of $C^{\infty}_{0}(\Omega)$ in $W$, that is, $$W_{0}=\overline{C^{\infty}_{0}(\Omega)}^{||.||_{W}} .$$
  \begin{enumerate}[label=(\roman*)]
  \item Theorem $\ref{1.1}$  remains true if we replace $W$ by $W_{0}$, $($$\Omega$ is not necessarily smooth$)$.
  \item  Since $ 1<p^{-}\leqslant \bar{p}(x)<p^{\ast}_{s}(x)$, for any $x\in \overline{\Omega}$, then Theorem $\ref{1.1}$  implies that $[.]_{s,p(x,y)}$ is a norm on
  $W_{0}$, which is equivalent to the norm  $\|.\|_W$. So $(W_{0},[.]_{s,p(x,y)})$ is a Banach space.
  
  \end{enumerate}
  \end{remark}
  Let denote by the $\mathcal{L}$ the operator associated to the $(-\Delta_{p(x,.)})^{s}$ defined as
  
  $$\mathcal{L}:W_{0}\longrightarrow W^{\ast}_{0}$$
  $$\hspace{3.3cm}u\longrightarrow\mathcal{L}(u) : W_{0}\longrightarrow \mathbb{R}$$
  $$\hspace{7.7cm}\varphi\longrightarrow\mathcal<\mathcal{L}(u),\varphi>$$
  such that $$<\mathcal{L}(u),\varphi>=\int_{\Omega\times\Omega}\frac{|u(x)-u(y)|^{p(x,y)-2}(u(x)-u(y))(\varphi(x)-\varphi(y))}{|x-y|^{N+sp(x,y)}}~dxdy,$$
  
  where $W^{\ast}_{0}$ is the dual space of $W_{0}$. \\

  \begin{lemma}\label{1.4}\textnormal{(\cite{2})}.
  Assume that hypothesis $(\ref{1})$ and $(\ref{2})$ are satisfied and $s\in (0,1)$. Then, the the following assertions hold:
  \begin{itemize}
  \item $\mathcal{L}$  is a bounded and strictly monotone operator.
  \item $\mathcal{L}$ is a mapping of type $(S_{+})$, that is,\\ if $u_{k}\rightharpoonup u$ in $W_{0}$ and $\underset{k\longrightarrow +\infty}{\limsup}<\mathcal{L}(u_{k})-\mathcal{L}(u), u_{k}-u>\leqslant0$, then  $u_{k}\longrightarrow u$ in $W_{0}$.
  \item $\mathcal{L}$ is a homeomorphism. 
  \end{itemize}
  \end{lemma}
  
\section{Functional framework}\label{300}
One of the aims of this paper is to study nonlocal problems driven by $\mathcal{L}^{p(x,.)}_K$ and  $(-\Delta_{p(x,.)})^{s}$  with Dirichlet boundary data via variational methods. For this purpose, we need to work in a suitable fractional Sobolev space. For this, we consider a functional analytical setting that is inspired by (but not equivalent to) the fractional Sobolev spaces in order to correctly encode the Dirichlet boundary datum in the variational formulation.\\ This section is devoted to the definition of this space as well as to its properties. Further,   we will prove a continuous compact embedding theorem of these spaces into variable exponent Lebesgue spaces. Finally, we establish a convergence property for a bounded sequence in $W^{K,p(x,y)}_{0}(\Omega)$.\\
Let $\Omega$ be a Lipschitz open bounded subset of $\mathbb{R}^N$,  $s\in(0,1)$ be fixed such that $sp^+<N$.  Denote
by $Q$ the set
$$Q:=\mathbb{R}^N\times\mathbb{R}^N\setminus (C\Omega\times C\Omega),\quad\text{where }~~ C\Omega=\mathbb{R}^N\setminus\Omega.$$

Now, due to the non-locality of the operator $\mathcal{L}^{p(x,.)}_K$  we  introduce the general fractional Sobolev space with variable exponent as follows 
$$W^{K,p(x,y)}(\Omega)=~ \left\{
   \begin{array}{clcclc}
   \hspace{-0.5cm}u:\mathbb{R}^N\longrightarrow\mathbb{R} ~\text{measurable, such that}~~u_{|\Omega}\in L^{\bar{p}(x)}(\Omega) ~\text{with}~\\   \int_{Q }\frac{|u(x)-u(y)|^{p(x,y)}}{\lambda^{p(x,y)}}K(x,y)~dxdy <+\infty,~~ \text{for some}~~\lambda>0\\
  
   \end{array}
   \right\}.
   $$
The norm in $W^{K,p(x,y)}(\Omega)$ can be defined as follows:
$$\|u\|_{W^{K,p(x,y)}(\Omega)}=\|u\|_{K,p(x,y)}=\|u\|_{L^{\bar{p}(x)}(\Omega)}+[u]_{K,p(x,y)},$$
where,  \small{$[u]_{K,p(x,y)}= \inf \bigg\{\lambda>0:\int_{Q}\frac{|u(x)-u(y)|^{p(x,y)}}{\lambda^{p(x,y)}}K(x,y)~dxdy \leqslant1 \bigg\}$}, (see Lemma $\ref{3.1}$).\\
For any $u\in W^{K,p(x,y)}(\Omega)$, we define the functional $$\rho_{K,p(.,.)}(u)=\int_{Q}|u(x)-u(y)|^{p(x,y)}K(x,y)~dxdy+\int_{\Omega}|u(x)|^{\bar{p}(x)}dx.$$
It is easy to see that $\rho_{K,p(.,.)}$ is a \textit{convex modular} on $ W^{K,p(x,y)}(\Omega)$. The norm associated with $\rho_{K,p(.,.)}$
 is given by $$\|u\|_{\rho_{K,p(.,.)}}=\inf \bigg\{\lambda>0 : \rho_{K,p(.,.)}\bigg(\frac{u}{\lambda}\bigg)\leqslant1\bigg\}.$$
Using the same argument as in \cite[Theorem 2.17]{7}, we prove that $\|.\|_{\rho_{K,p(.,.)}}$ is a norm on $W^{K,p(x,y)}(\Omega)$, which is equivalent to the norm $\|.\|_{K,p(x,y)}$. 
We also define the closed linear subspace of $W^{K,p(x,y)}(\Omega)$ by 
$$W_{0}^{K,p(x,y)}(\Omega)=\bigg\{u\in W^{K,p(x,y)}(\Omega)~ :~ u(x)=0~~ \text{a.e.} ~~x\in \mathbb{R}^{N}\setminus \Omega\bigg\}.$$
On the other hand, for any $u\in W_{0}^{K,p(x,y)}(\Omega)$, we define the functional $$\rho^{o}_{K,p(.,.)}(u)=\int_{Q}|u(x)-u(y)|^{p(x,y)}K(x,y)~dxdy.$$
$\rho^{o}_{K,p(.,.)}$ is a \textit{convex modular} on $ W_{0}^{K,p(x,y)}(\Omega)$. The norm associated with $\rho^{o}_{K,p(.,.)}$
 is given by $$\|u\|_{\rho^{o}_{K,p(.,.)}}=[u]_{K,p(x,y)}=inf \bigg\{\lambda>0 : \rho^{o}_{K,p(.,.)}\bigg(\frac{u}{\lambda}\bigg)\leqslant1\bigg\}.$$

\begin{remark}\label{3.4}\textbf{}
\begin{enumerate}
\item  $\rho^{o}_{K,p(.,.)}$  also check the results of Proposition $\ref{2.2.2}$.
\item The modular  $\rho^{o}_{K,p(.,.)}$ does not satisfy the triangle inequality, that is,$$\rho^{o}_{K,p(.,.)}(u+v)\leqslant \rho^{o}_{K,p(.,.)}(u)+\rho^{o}_{K,p(.,.)}(v)$$
However, there is a substitute that is sometimes useful. 
$$\rho^{o}_{K,p(.,.)}(u+v)\leqslant 2^{p^{+}-1} \big(\rho^{o}_{K,p(.,.)}(u)+\rho^{o}_{K,p(.,.)}(v)\big).$$
We will refer to this as the modular triangle inequality. 
\end{enumerate}
\end{remark}
\begin{lemma}\label{3.1}
$\|.\|_{K,p(x,y)}$ is a norm on $W^{K,p(x,y)}(\Omega)$.  
\end{lemma}
\textbf{\textit{Proof}}. Since $\|.\|_{L^{\bar{p}(x)}(\Omega)}$ is a norm on $L^{\bar{p}(x)}(\Omega)$. So we need to prove that: \\ 
$(i)$ $\|u\|_{K,p(x,y)}=0$ if and only if $u=0$,\\
and $[.]_{K,p(x,y)}$ is a semi-norm on $W^{K,p(x,y)}(\Omega)$, that is, \\
$(ii)$ for  all $\alpha\in \mathbb{R}$, $[\alpha u]_{K,p(x,y)}= |\alpha|[u]_{K,p(x,y)}$,\\
$(iii)$ $[u+v]_{K,p(x,y)}\leqslant [u]_{K,p(x,y)}+[v]_{K,p(x,y)}$.\\

Indeed, For $(i)$, if $u=0$, then  $$\int_{Q}|u(x)-u(y)|^{p(x,y)}K(x,y)~dxdy=0 \hspace{0.3cm}and\hspace{0.3cm} \|u\|_{L^{\bar{p}(x)}(\Omega)}=0.$$
So, by Lemma \ref{2.3} we get  $$[u]_{K,p(x,y)}=0\hspace{0.3cm}\text{and }\hspace{0.3cm}\|u\|_{L^{\bar{p}(x)}(\Omega)}=0\hspace{0.3cm}\Rightarrow\hspace{0.3cm}\|u\|_{K,p(x,y)}=0.$$

Conversely, if $\|u\|_{K,p(x,y)}=0$, then
\begin{equation}\label{8}
\|u\|_{L^{\bar{p}(x)}(\Omega)}=0,
\end{equation}
and
\begin{equation}\label{9}
[u]_{K,p(x,y)}=0.
\end{equation}
By (\ref{8}), we have 
\begin{equation}\label{10}
u=0 \hspace{0.3cm}\text{a.e. in}\hspace{0.3cm} \Omega.
\end{equation}
Using (\ref{9}), we get $$\lambda^{\star}=\inf \bigg\{\lambda>0:\int_{Q}\frac{|u(x)-u(y)|^{p(x,y)}}{\lambda^{p(x,y)}}K(x,y)~dxdy \leqslant1 \bigg\}=0.$$
Let $\lambda_n>0$, $n\in \mathbb{N}$, such that $\lambda_n$ decreases to $\lambda^{\star}$ and $$\int_{Q}\frac{|u(x)-u(y)|^{p(x,y)}}{\lambda_{n}^{p(x,y)}}K(x,y)~dxdy \leqslant1.$$
If $\lambda_{n}<1$, then wa have $$\int_{Q}\frac{|u(x)-u(y)|^{p(x,y)}}{\lambda_{n}^{p^{-}}}K(x,y)~dxdy \leqslant1.$$
Hence $$\int_{Q}|u(x)-u(y)|^{p(x,y)}K(x,y)~dxdy \leqslant\lambda_{n}^{p^{-}}.$$
When $n\longrightarrow +\infty$, we obtain $\lambda_{n}^{p^{-}}\longrightarrow\lambda^{\star}=0$. Thus
$$0\leqslant\int_{Q}|u(x)-u(y)|^{p(x,y)}K(x,y)~dxdy\leqslant0\hspace{0.2cm}\Rightarrow\int_{Q}|u(x)-u(y)|^{p(x,y)}K(x,y)~dxdy=0.$$
We conclude that $u(x)=u(y)$ a.e. $(x,y)\in Q$, then $u=c\in\mathbb{R}$ a.e. in $\mathbb{R}^N$.
Finally, by (\ref{10}) it easily follows that $c=0$, so $u=0$ a.e. in $\mathbb{R}^N$.\\
To prove $(ii)$, note that if $\alpha=0$, this follows from $(i)$. Fix $\alpha\neq0$, then by a change of variable, we have 
$$\begin{array}{clll}
[\alpha u]_{K,p(x,y)}&=&
\inf \bigg\{\lambda>0:\int_{Q}\bigg|\frac{ \alpha u(x)-\alpha u(y)}{\lambda}\bigg|^{p(x,y)}K(x,y)dxdy\leqslant1\bigg\} \\\\&=& \inf \bigg\{\lambda>0:\int_{Q}\frac{| u(x)-u(y)|^{p(x,y)}}{\big(\frac{\lambda}{|\alpha|}\big)^{p(x,y)}}K(x,y)dxdy\leqslant1\bigg\}
\\\\&=& |\alpha|\inf \bigg\{\frac{\lambda}{|\alpha|}>0:\int_{Q}\frac{| u(x)-u(y)|^{p(x,y)}}{\big(\frac{\lambda}{|\alpha|}\big)^{p(x,y)}}K(x,y)dxdy\leqslant1\bigg\}
\\\\ &=& |\alpha|\inf \bigg\{\mu>0:\int_{Q}\frac{| u(x)-u(y)|^{p(x,y)}}{\mu^{p(x,y)}}K(x,y)dxdy\leqslant1\bigg\}\\\\&=&|\alpha| [ u]_{K,p(x,y)}.
\end{array}
$$

Finally, to prove $(iii)$, fix $\lambda_{u}>[ u]_{K,p(x,y)}$ and $\lambda_{v}>[ v]_{K,p(x,y)}$. Then  $$\rho^{o}_{K,p(.,.)}\bigg(\frac{u}{\lambda_{u}}\bigg)\leqslant1 \hspace{0.5cm}and\hspace{0.5cm} \rho^{o}_{K,p(.,.)}\bigg(\frac{v}{\lambda_{v}}\bigg)\leqslant1$$
Now, Let $\lambda=\lambda_{u}+\lambda_{v}$, then by the convexity  of $\rho^{o}_{K,p(.,.)}$, we have 
$$\begin{array}{clll}
\rho^{o}_{K,p(.,.)}\bigg(\frac{u+v}{\lambda}\bigg)&=& \rho^{o}_{K,p(.,.)}\bigg(\frac{\lambda_{u}u}{\lambda\lambda_{u}}+\frac{\lambda_{v}v}{\lambda\lambda_{v}}\bigg)\\\\&\leqslant&\frac{\lambda_{u}}{\lambda}\rho^{o}_{K,p(.,.)}\bigg(\frac{u}{\lambda_{u}}\bigg)+ \frac{\lambda_{v}}{\lambda}\rho^{o}_{K,p(.,.)}\bigg(\frac{v}{\lambda_{v}}\bigg)
\end{array}$$
Since, $\frac{\lambda_{u}}{\lambda}+\frac{\lambda_{v}}{\lambda}=1$, then $$\rho^{o}_{K,p(.,.)}\bigg(\frac{u+v}{\lambda}\bigg)\leqslant1.$$
Hence, $$[u+v]_{K,p(x,y)}\leqslant\lambda=\lambda_{u}+\lambda_{v},$$
we take the  infimum over all such $\lambda_{u}$ and
$\lambda_{v}$, we get the desired inequality. \hspace{2cm}$\Box$
 
\begin{remark}
We remark that in the model case in which $K(x,y)= |x-y|^{-(N+sp(x,y))}$ the norms $\|.\|_{K,p(x,y)}$ and $\|.\|_{s,p(x,y)}$ are not the same, because $\Omega\times\Omega$ is strictly contained in $Q$: This makes the fractional Sobolev space with variable exponent $W^{s,p(x,y)}(\Omega)$ not sufficient for studying the nonlocal problems.
\end{remark}

\begin{lemma}\label{3.2}
Let $K:\mathbb{R}^N \times\mathbb{R}^N\longrightarrow (0,+\infty)$ be a measurable function satisfy $(\ref{4})$ and $(\ref{6})$, let  $p:\mathbb{R}^N \times\mathbb{R}^N\longrightarrow(1,+\infty)$ be a continuous bounded function satisfy $(\ref{1})$ and  $(\ref{2})$. Then
$$C^{\infty}_{0}(\Omega)\subset W_{0}^{K,p(x,y)}(\Omega).$$
\end{lemma}
\textbf{\textit{Proof}}. Using the same argument as in \cite{19}, this lemma can be proved. For completeness, we give its proof. For $u \in C^{\infty}_{0}(\Omega)$, we only need to check that $[u]_{K,p(x,y)}<+\infty$.
From Lemma \ref{2.3}, we need to prove that $$\int_{\mathbb{R}^{2N}}|u(x)-u(y)|^{p(x,y)}K(x,y)~dxdy <+\infty.$$
Indeed, since $u=0$ in $\mathbb{R}^N \setminus \Omega$, we have that
\begin{equation}\label{11}
\begin{array}{clll}
\int_{\mathbb{R}^{2N}}|u(x)-u(y)|^{p(x,y)}K(x,y)~dxdy &=&\hspace{-0.8cm}\int_{\Omega\times\Omega}|u(x)-u(y)|^{p(x,y)}K(x,y)~dxdy  
\\\\&\hspace{0.6cm}+& 2\int_{\Omega\times(\mathbb{R}^N \setminus \Omega)}|u(x)-u(y)|^{p(x,y)}K(x,y)~dxdy \\&\leqslant&\hspace{-0.8cm}2\int_{\Omega\times\mathbb{R}^{N}}|u(x)-u(y)|^{p(x,y)}K(x,y)~dxdy.
\end{array}
\end{equation}
Notice that $$|u(x)-u(y)|\leqslant \|\nabla u\|_{L^{\infty}(\mathbb{R}^{N})} |x-y|\hspace{0.2cm}\text{and}\hspace{0.2cm}|u(x)-u(y)|\leqslant 2\| u\|_{L^{\infty}(\mathbb{R}^{N})},\hspace{0.1cm}\text{ for all } (x,y)\in \mathbb{R}^{2N}.$$
Hence,
$$|u(x)-u(y)|^{p(x,y)}\leqslant \|\nabla u\|^{p(x,y)}_{L^{\infty}(\mathbb{R}^{N})} |x-y|^{p(x,y)},\text{ for all } (x,y)\in \mathbb{R}^{2N}$$and$$|u(x)-u(y)|^{p(x,y)}\leqslant 2^{p(x,y)}\| u\|^{p(x,y)}_{L^{\infty}(\mathbb{R}^{N})},\hspace{0.1cm}\text{ for all } (x,y)\in \mathbb{R}^{2N}.$$
Thus, 
$$|u(x)-u(y)|^{p(x,y)}\leqslant \bigg( \|\nabla u\|^{p^{+}}_{L^{\infty}(\mathbb{R}^{N})}+\|\nabla u\|^{p^{-}}_{L^{\infty}(\mathbb{R}^{N})}\bigg) |x-y|^{p(x,y)},\text{ for all } (x,y)\in \mathbb{R}^{2N}$$and$$|u(x)-u(y)|^{p(x,y)}\leqslant 2^{p+}\bigg(\| u\|^{p^{+}}_{L^{\infty}(\mathbb{R}^{N})}+\| u\|^{p^{-}}_{L^{\infty}(\mathbb{R}^{N})}\bigg),\hspace{0.1cm}\text{ for all } (x,y)\in \mathbb{R}^{2N}.$$
So, we have 
$$|u(x)-u(y)|^{p(x,y)}\leqslant 2^{p^{+}+1}\bigg(\| u\|^{p^{+}}_{C^{1}(\mathbb{R}^{N})}+\| u\|^{p^{-}}_{C^{1}(\mathbb{R}^{N})}\bigg)\hspace{0.1cm}\min\bigg\{1,|x-y|^{p(x,y)}\bigg\}.$$
From (\ref{11}), we deduce that 
\begingroup\makeatletter\def\f@size{9.3}\check@mathfonts$$\begin{array}{clll}
\int_{\mathbb{R}^{2N}}|u(x)-u(y)|^{p(x,y)}K(x,y)~dxdy &=&2^{p^{+}+1}\bigg(\| u\|^{p^{+}}_{C^{1}(\mathbb{R}^{N})}+\| u\|^{p^{-}}_{C^{1}(\mathbb{R}^{N})}\bigg)\int_{\Omega\times\mathbb{R}^{N}}m(x,y)K(x,y)~dxdy
\\\\&\leqslant&C_{1}\int_{\mathbb{R}^{N}\times\mathbb{R}^{N}}m(x,y)K(x,y)~dxdy
\end{array}$$\endgroup
Assumption (\ref{6}) implies that $u\in W_{0}^{K,p(x,y)}(\Omega).$ \hspace{6.5cm}$\Box$
\begin{remark}
A trivial consequence of Lemma $\ref{3.2}$, $W^{K,p(x,y)}(\Omega)$ and $W_{0}^{K,p(x,y)}(\Omega)$ are non-empty.
\end{remark}
The modular $\rho^{o}_{K,p(.,.)}$  check the following result, which is similar to Proposition \ref{2.2.1} and Lemma \ref{2.2}.
   \begin{lemma}\label{2.3}
    Let $p:\mathbb{R}^N\times\mathbb{R}^N\longrightarrow(1,+\infty)$
     be a continuous variable exponent and  $K:\mathbb{R}^N \times\mathbb{R}^N\longrightarrow (0,+\infty)$ is a measurable function satisfy $(\ref{4})$ and $(\ref{6})$. Then
    For any $u\in W^{K,p(x,y)}_0$, we have
    \begin{enumerate}[label=(\roman*)]
    \small \item $1\leqslant [u]_{K,p(x,y)}$ $\Rightarrow$ $[u]^{p^{-}}_{K,p(x,y)}\leqslant\rho^{o}_{K,p(.,.)}(u)\leqslant [u]^{p^{+}}_{K,p(x,y)}$,
    \item  $ [u]_{K,p(x,y)}\leqslant1$ $\Rightarrow$ $[u]^{p^{+}}_{K,p(x,y)}\leqslant\rho^{o}_{K,p(.,.)}(u)\leqslant [u]^{p^{-}}_{K,p(x,y)}$.
    \end{enumerate}
    \end{lemma} 
    
\textbf{\textit{Proof}}. We prove the first pair of inequalities; the proof of the second is essentially
the same. Indeed, it is easy to see that, for all 
 $\lambda\in (0,1)$, we get  
 $$\lambda^{p^{+}}\rho^{o}_{K,p(.,.)}(u)\leqslant \rho^{o}_{K,p(.,.)}(\lambda u)\leqslant \lambda^{p^{-}}\rho^{o}_{K,p(.,.)}(u)$$
 Now, if  $[u]_{K,p(x,y)}>1$, then $0<\frac{1}{[u]_{K,p(x,y)}}<1$, so we have 
 $$\frac{\rho^{o}_{K,p(.,.)}(u)}{[u]^{p^{+}}_{K,p(x,y)}}\leqslant\rho^{o}_{K,p(.,.)}\bigg(\frac{u}{[u]_{K,p(x,y)}}\bigg)\leqslant\frac{\rho^{o}_{K,p(.,.)}(u)}{[u]^{p^{-}}_{K,p(x,y)}}.$$
 Since $\rho^{o}_{K,p(.,.)}\bigg(\frac{u}{[u]_{K,p(x,y)}}\bigg)=1$, so the desired result follows. \hspace{4.5cm}$\Box$
 
In the following lemma we compare the spaces $W^{K,p(x,y)}$ and $W^{s,p(x,y)}$. This lemma is  crucial in the proof of  the continuous embedding theorem.
\begin{lemma}\label{3.3}
Let $K:\mathbb{R}^N \times\mathbb{R}^N\longrightarrow (0,+\infty)$ be a measurable function satisfying $(\ref{4})$-$(\ref{6})$. Let $p:\mathbb{R}^N \times\mathbb{R}^N\longrightarrow(1,+\infty)$ be a continuous bounded function satisfying $(\ref{1})$ and $(\ref{2})$. Then the following assertions hold:
\begin{enumerate}[label=(\roman*)]
\item If $u\in W^{K,p(x,y)}(\Omega)$, then $u\in W^{s,p(x,y)}(\Omega)$. Moreover, $$\|u\|_{s,p(x,y)}\leqslant \max\big\{1, \tilde{k}_0\big\}\|u\|_{K,p(x,y)},$$
where $\tilde{k}_0=\tilde{k}_0(k_0, p^{-},  p^{+})$ is a positive constant.  That is, the space $W^{K,p(x,y)}(\Omega)$ is continuously embedded in $W^{s,p(x,y)}(\Omega)$.
\item If $u\in W_{0}^{K,p(x,y)}(\Omega)$, then $u\in W^{s,p(x,y)}(\mathbb{R}^{N})$. Moreover, $$\|u\|_{W^{s,p(x,y)}(\Omega)}\leqslant\|u\|_{W^{s,p(x,y)}(\mathbb{R}^N)}\leqslant \max\big\{1, \tilde{k}_0\big\}\|u\|_{K,p(x,y)}$$
\end{enumerate}
\end{lemma}
\textbf{\textit{Proof}}. $(i)$- Let $\lambda>0$, for $u\in W^{K,p(x,y)}(\Omega)$ and by (\ref{5}), we have 
\begin{equation}\label{12}
\int_{\Omega\times\Omega}\frac{|u(x)-u(y)|^{p(x,y)}}{\lambda^{p(x,y)}|x-y|^{sp(x,y)+N}}~dxdy \leqslant\int_{Q}\frac{|u(x)-u(y)|^{p(x,y)}}{\lambda^{p(x,y)}}\frac{K(x,y)}{k_0}~dxdy
.\end{equation}
We define 
$$\mathcal{A}^{s}_{\lambda,\Omega}=\bigg\{\lambda>0:\int_{\Omega\times\Omega}\frac{|u(x)-u(y)|^{p(x,y)}}{\lambda^{p(x,y)}|x-y|^{sp(x,y)+N}}~dxdy\leqslant1 \bigg\},$$
and
$$\mathcal{A}^{K,k_0}_{\lambda,Q}=\bigg\{\lambda>0:\int_{Q}\frac{|u(x)-u(y)|^{p(x,y)}}{\lambda^{p(x,y)}}\frac{K(x,y)}{k_0}~dxdy\leqslant1 \bigg\}.$$
By (\ref{12}), it is easy to see that $\mathcal{A}^{K,k_0}_{\lambda,Q}\subset\mathcal{A}^{s}_{\lambda,\Omega}$. Hence $\inf_{\lambda>0}\mathcal{A}^{s}_{\lambda,\Omega}\leqslant\inf_{\lambda>0}\mathcal{A}^{K,k_0}_{\lambda,Q}$. Then, we have
\begin{equation}\label{13}
[u]_{s,p(x,y)}=\inf_{\lambda>0}\mathcal{A}^{s}_{\lambda,\Omega}\leqslant\inf_{\lambda>0}\mathcal{A}^{K,k_0}_{\lambda,Q}.
\end{equation}
Now, let $$\mathcal{A}^{K}_{\lambda,Q}=\bigg\{\lambda>0:\int_{Q}\frac{|u(x)-u(y)|^{p(x,y)}}{\lambda^{p(x,y)}}K(x,y)~dxdy\leqslant1 \bigg\}.$$
And we set $$\tilde{k}_0=\max\bigg\{k_{0}^{-\frac{1}{p^{-}}}, k_{0}^{-\frac{1}{p^{+}}}\bigg\}\hspace{0.4cm}and \hspace{0.4cm} \bar{\lambda}= \lambda \tilde{k}_{0}$$
So ,  we obtain \begin{equation}
\begin{array}{clll}
\tilde{k}_0 \inf_{\lambda>0}\mathcal{A}^{K}_{\lambda,Q}&=&\inf_{\lambda>0}
\bigg\{\lambda\tilde{k}_0:\int_{Q}\frac{|u(x)-u(y)|^{p(x,y)}}{\lambda^{p(x,y)}}K(x,y)~dxdy\leqslant1 \bigg\}\\\\&=&\inf_{\bar{\lambda}>0}
\bigg\{\bar{\lambda}:\int_{Q}\frac{|u(x)-u(y)|^{p(x,y)}}{\bar{\lambda}^{p(x,y)}}\tilde{k}_{0}^{p(x,y)}K(x,y)~dxdy\leqslant1 \bigg\}
\end{array}
\end{equation}\label{14}
Since $\tilde{k}_{0}\geqslant k_{0}^{-\frac{1}{p(x,y)}}$, then $\tilde{k}_{0}^{p(x,y)}\geqslant \frac{1}{k_{0}}$. So, we get
$$\int_{Q}\frac{|u(x)-u(y)|^{p(x,y)}}{\lambda^{p(x,y)}}\frac{K(x,y)}{k_{0}}~dxdy \leqslant\int_{Q}\frac{|u(x)-u(y)|^{p(x,y)}}{\lambda^{p(x,y)}}\tilde{k}_{0}^{p(x,y)}K(x,y)~dxdy$$ 
Let 
$$\mathcal{B}^{K,\tilde{k}_{0}}_{\lambda,Q}=\bigg\{\lambda>0:\int_{Q}\frac{|u(x)-u(y)|^{p(x,y)}}{\lambda^{p(x,y)}}\tilde{k}_{0}^{p(x,y)}K(x,y)~dxdy\leqslant1 \bigg\}.$$
We remark that $\mathcal{B}^{K,\tilde{k}_{0}}_{\lambda,Q}\subset\mathcal{A}^{K,k_0}_{\lambda,Q}$. This implies that $\inf_{\lambda>0}\mathcal{A}^{K,k_0}_{\lambda,Q}\leqslant\inf_{\lambda>0}\mathcal{B}^{K,\tilde{k}_{0}}_{\lambda,Q}$.\\ 
Using \eqref{14}, we get 
$$\inf_{\lambda>0}\mathcal{A}^{K,k_0}_{\lambda,Q}\leqslant\inf_{\lambda>0}\mathcal{B}^{K,\tilde{k}_{0}}_{\lambda,Q}=\tilde{k}_0\inf_{\lambda>0}\mathcal{A}^{K}_{\lambda,Q} $$
Hence,by (\ref{13}), we have
$$[u]_{s,p(x,y)}\leqslant\tilde{k}_0[u]_{K,p(x,y)}<+\infty.$$
In fact, by using the definition of norms $\|u\|_{s,p(x,y)}$ and  $\|u\|_{K,p(x,y)}$, we infer  $$\|u\|_{s,p(x,y)}\leqslant max\big\{1, \tilde{k}_0\big\}\|u\|_{K,p(x,y)}.$$
The first assertion is proved.\\
$(ii)$-  For $u\in W_0^{K,p(x,y)}(\Omega)$, we have $u=0$ a.e. in $\mathbb{R}^{N}\setminus\Omega$. Then 
$$\|u\|_{L^{\bar{p}(x)}(\Omega)}=\|u\|_{L^{\bar{p}(x)}(\mathbb{R}^{N})}<+\infty.$$
 By the same argument in the assertion $(i)$, we get $u\in W^{s,p(x,y)}(\mathbb{R}^{N})$ and $$[u]_{W^{s,p(x,y)}(\Omega)}\leqslant[u]_{W^{s,p(x,y)}(\mathbb{R}^{N})}\leqslant\tilde{k}_0[u]_{K,p(x,y)}.$$
So the estimate on the norm is easily follows.\hspace{7cm}$\Box$\\

Now, we are ready to prove the main theorem of this section.
\begin{theorem}\label{3.1.1}
Let $\Omega$ be a Lipschitz  bounded domain in $\mathbb{R}^{N}$ and $s\in (0,1)$. Let $p:\mathbb{R}^{N}\times\mathbb{R}^{N}\longrightarrow(1,+\infty)$
 be a continuous variable exponent satisfies $(\ref{1})$ and  $(\ref{2})$ with $s p^+<N$.  Let $r:\overline{\Omega}\longrightarrow(1,+\infty)$  be a continuous bounded variable exponent such that
  $$1<r^-\leqslant r(x)<p^{\ast}_{s}(x), ~~ \text{ for all } x\in\overline{\Omega}.$$
  Suppose that  $K:\mathbb{R}^N \times\mathbb{R}^N\longrightarrow (0,+\infty)$ is a measurable function satisfying $(\ref{4})$-$(\ref{6})$. Then
  \begin{enumerate}[label=(\roman*)]
  \item There exists a positive constant $C=C(N,p,r,s,\Omega)>0$, such that for any $u\in W^{K,p(x,y)}(\Omega)$, we have
  $$\|u\|_{L^{r(x)}(\Omega)}\leqslant C\|u\|_{W^{s,p(x,y)}(\Omega)}\leqslant C\max\big\{1, \tilde{k}_0\big\}\|u\|_{K,p(x,y)},$$
   That is, the space $W^{K,p(x,y)}(\Omega)$ is continuously embedded in $L^{r(x)}(\Omega)$. Moreover, this
    embedding is compact.
  \item There exists a positive constant $C_0=C_0(N,p,s,\tilde{k}_0,\Omega)>0$, such that$$[u]_{K,p(x,y)}\leqslant\|u\|_{K,p(x,y)}\leqslant C_0[u]_{K,p(x,y)}.$$
  \end{enumerate}
\end{theorem}
\textbf{\textit{Proof}}. $(i)$- Let $u\in W^{K,p(x,y)}(\Omega)$, by Lemme \ref{3.3}, we have $u\in W^{s,p(x,y)}(\Omega)$ and
\begin{equation}\label{15}
\|u\|_{W^{s,p(x,y)}}=\|u\|_{s,p(x,y)}\leqslant max\big\{1, \tilde{k}_0\big\}\|u\|_{K,p(x,y)}.
\end{equation}
Combining (\ref{15}) with Theorem \ref{1.1}, we obtain 
 $$\|u\|_{L^{r(x)}(\Omega)}\leqslant C\|u\|_{W^{s,p(x,y)}(\Omega)}\leqslant C\max\big\{1, \tilde{k}_0\big\}\|u\|_{K,p(x,y)}.$$
 Hence, we deduce that
 $$W^{K,p(x,y)}(\Omega)\hookrightarrow W^{s,p(x,y)}(\Omega)\hookrightarrow \hookrightarrow L^{r(x)}(\Omega).$$
 Since the latter embedding is compact, then the embedding $W^{K,p(x,y)}(\Omega)\hookrightarrow L^{r(x)}(\Omega)$ is also compact.\\
 
 $(ii)$- This assertion is easily follows by combining the definition of $\|.\|_{K,p(x,y)}$ with assertion $(i)$ and assumptions (\ref{4})-(\ref{6}).\hspace{8.7cm}$\Box$
 \begin{remark}\label{r3.4}\textit{}
\begin{enumerate}
\item The assertion $(i)$ implies also that   $W^{K,p(x,y)}_0(\Omega)$ is  continuously embedded in $L^{r(x)}(\Omega)$, where  $1<r^-\leqslant r(x)<p^{\ast}_{s}(x)$ for any $ x\in\overline{\Omega}$. Moreover, this
    embedding is compact.
\item As a consequence of assertion $(ii)$, $[.]_{K,p(x,y)}$ is an equivalent norm of $\|u\|_{K,p(x,y)}$  on $W_0^{K,p(x,y)}(\Omega)$.

\end{enumerate}
 \end{remark}
 \begin{lemma}\label{3.4.1}
 $\left( W^{K,p(x,y)}_0(\Omega), [.]_{K,p(x,y)}\right) $ is a separable, reflexive, and uniformly convex Banach space.
 \end{lemma}
 \textbf{\textit{Proof}}.
 We first prove that $W^{K,p(x,y)}_0(\Omega)$ is complete with respect to the norm $[.]_{K,p(x,y)}$.
 Let $\{u_n\}$ be a Cauchy sequence in $W^{K,p(x,y)}_0(\Omega)$.  Since $\bar{p}(x)<p^{\ast}_{s}(x)$, so, combining ($i$) and ($ii$) of Theorem \ref{3.1.1}, for any $\varepsilon>0$, there exists $n^{\star}_{\varepsilon}$ such that if $n,m\geqslant n^{\star}_{\varepsilon}$, we get 
 \begin{equation}\label{16}
\frac{1}{\overline{C}} \|u_n-u_m\|_{L^{\bar{p}(x)}(\Omega)}\leqslant [u_n-u_m]_{K,p(x,y)}\leqslant\varepsilon,
 \end{equation}
 where $\overline{C}= C_0C  max\big\{1, \tilde{k}_0\big\}$. By the completeness of  $L^{\bar{p}(x)}(\Omega)$, there exists $u\in L^{\bar{p}(x)}(\Omega)$ such that $u_n\longrightarrow u$ strongly in $L^{\bar{p}(x)}(\Omega)$ as $n\longrightarrow+\infty$. Since $u_n=0$ a.e. in $\mathbb{R}^{N}\setminus\Omega$, so we define $u=0$  a.e. in $\mathbb{R}^{N}\setminus\Omega$. Then $u_n\longrightarrow u$ strongly in $L^{\bar{p}(x)}(\mathbb{R}^{N})$ as $n\longrightarrow+\infty$. So there exists a subsequence $\{u_{n_{j}}\}$ of $\{u_n\}$ in $W^{K,p(x,y)}_0(\Omega)$, such that $u_{n_{j}}\longrightarrow u$  a.e. in $\mathbb{R}^{N}$.\\
 Now, we need to prove that $[u]_{K,p(x,y)}<+\infty$. By Lemma \ref{2.3}, is enough to show that $$\rho^{o}_{K,p(.,.)}(u)<+\infty.$$
 Indeed, by the Fatou Lemma, with $\varepsilon=1$, we have 
 $$\begin{array}{clll}
 \rho^{o}_{K,p(.,.)}(u)&=&\hspace{-6cm}\int_{Q}|u(x)-u(y)|^{p(x,y)}K(x,y)~dxdy 
  \\\\&\leqslant&\hspace{-6cm}\liminf_{j\longrightarrow+\infty}\int_{Q}|u_{n_{j}}(x)-u_{n_{j}}(y)|^{p(x,y)}K(x,y)~dxdy
  \\\\&\leqslant&\hspace{-6cm}\liminf_{j\longrightarrow+\infty}\int_{Q}|u_{n_{j}}(x)-u_{n_{j}}(y)-\big(u_{n^{\star}_{1}}(x)-u_{n^{\star}_{1}}(y)\big)+\big(u_{n^{\star}_{1}}(x)-u_{n^{\star}_{1}}(y)\big)|^{p(x,y)}K(x,y)~dxdy
  \\\\&\leqslant&\hspace{-6cm}\liminf_{j\longrightarrow+\infty}\int_{Q}2^{p^{+}-1}\bigg(\big|u_{n_{j}}(x)-u_{n_{j}}(y)-\big(u_{n^{\star}_{1}}(x)-u_{n^{\star}_{1}}(y)\big)\big|^{p(x,y)}\\&\hspace{6cm}+&\big|\big(u_{n^{\star}_{1}}(x)-u_{n^{\star}_{1}}(y)\big)\big|^{p(x,y)}\bigg)K(x,y)~dxdy
 \\\\&\leqslant&\hspace{-6cm}2^{p^{+}-1}\liminf_{j\longrightarrow+\infty}\bigg\{\int_{Q}|(u_{n_{j}}-u_{n^{\star}_{1}})(x)-(u_{n_{j}}-u_{n^{\star}_{1}})(y)|^{p(x,y)}K(x,y)~dxdy\\&\hspace{5cm}+&\hspace{-1.2cm}\int_{Q}|u_{n^{\star}_{1}}(x)-u_{n^{\star}_{1}}(y)\big|^{p(x,y)}K(x,y)~dxdy \bigg\}
 \\&\leqslant&\hspace{-6cm} 2^{p^{+}-1}\bigg\{\liminf_{j\longrightarrow+\infty}\rho^{o}_{K,p(.,.)}(u_{n_{j}}-u_{n^{\star}_{1}})+\rho^{o}_{K,p(.,.)}(u_{n^{\star}_{1}})\bigg\}.
  
 \end{array}$$
 Using Lemma \ref{2.3} we obtain
 $$\rho^{o}_{K,p(.,.)}(u)\leqslant 2^{p^{+}-1}\bigg\{\liminf_{j\longrightarrow+\infty}\bigg([u_{n_{j}}-u_{n^{\star}_{1}}]^{p^{+}}_{K,p(x,y)}+[u_{n_{j}}-u_{n^{\star}_{1}}]^{p^{-}}_{K,p(x,y)}\bigg)+\bigg([u_{n^{\star}_{1}}]^{p^{+}}_{K,p(x,y)}+[u_{n^{\star}_{1}}]^{p^{-}}_{K,p(x,y)}\bigg)\bigg\}.$$
 By the inequality (\ref{16}) with $\varepsilon=1$, we get
 
  $$\rho^{o}_{K,p(.,.)}(u)\leqslant 2^{p^{+}-1}\bigg(2+[u_{n^{\star}_{1}}]^{p^{+}}_{K,p(x,y)}+[u_{n^{\star}_{1}}]^{p^{-}}_{K,p(x,y)}\bigg)<+\infty.$$
  Thus $u\in W^{K,p(x,y)}_0(\Omega)$. \\
  On the other hand, let $n\geqslant n^{\star}_{\varepsilon}$, combining (\ref{16}) with Lemma \ref{2.3} and the Fatou Lemma, we have
 $$\rho^{o}_{K,p(.,.)}(u_{n}-u)\leqslant\liminf_{j\longrightarrow+\infty}\rho^{o}_{K,p(.,.)}(u_n-u_{n_{j}})\leqslant \frac{\varepsilon^{p^{+}}+\varepsilon^{p{-}}}{2}=\varepsilon'.$$
 Hence, $$\lim_{n\longrightarrow+\infty}\rho^{o}_{K,p(.,.)}(u_{n}-u)=0.$$
 Using Remark \ref{3.4}-(1), we conclude that
  $$\lim_{n\longrightarrow+\infty}[u_{n}-u]_{K,p(x,y)}=0.$$
  That is, $u_n\longrightarrow u$ strongly in $W^{K,p(x,y)}_0(\Omega)$, as $n\longrightarrow+\infty$.\\
  Let us now prove that the space $W^{K,p(x,y)}_0(\Omega)$ is a separable and uniformly convex reflexive space. For this, we define the operator
  $$\mathcal{P}:W^{K,p(x,y)}_0(\Omega)\longrightarrow L^{p(x,y)}(Q,dxdy)$$
 $$\hspace{4.3cm}u\longrightarrow\big(u(x)-u(y)\big)K(x,y)^{\frac{1}{p(x,y)}}$$ 
 Clearly $\mathcal{P}$ is an isometry from $W^{K,p(x,y)}_0(\Omega)$ into $L^{p(x,y)}(Q)$. Since $W^{K,p(x,y)}_0(\Omega)$ is a Banach space, then $\mathcal{P}(W^{K,p(x,y)}_0(\Omega))$ is a closed subset of $L^{p(x,y)}(Q)$ (which is a separable and reflexive uniformly convex space, see Proposition \ref{2.3.1}). It follows that $\mathcal{P}(W^{K,p(x,y)}_0(\Omega))$ is separable and reflexive uniformly convex space.
  Consequently, $W^{K,p(x,y)}_0(\Omega)$ is also a separable and reflexive uniformly convex space.\\
  This concludes the proof.\hspace{10.5cm}$\Box$
\begin{corollary}\textbf{}
 \begin{enumerate}[label=(\roman*)]
 \item $\left( W^{K,p(x,y)}(\Omega), \|.\|_{K,p(x,y)}\right) $ is a separable and reflexive uniformly convex space.
 \item If $\Omega\subset\mathbb{R}^{N}$ is a domain of class $C^{0,1}$, then $\big(W^{K,p(x,y)}(\Omega), \|.\|_{K,p(x,y)}\big)$ is a Banach space. 
 \end{enumerate}
\end{corollary}
  \textbf{\textit{Proof}}.
  $(i)$- We consider the operator 
   $$\widetilde{\mathcal{P}}:W^{K,p(x,y)}(\Omega)\longrightarrow L^{\bar{p}}(\Omega)\times L^{p(x,y)}(Q,dxdy)=E$$
   $$\hspace{2.9cm}u\longrightarrow\big(u(x),(u(x)-u(y))K(x,y)^{\frac{1}{p(x,y)}}\big),$$ 
   which is an isometry from  $W^{K,p(x,y)}(\Omega)$ to $E$. The rest of proof is similar to Lemma \ref{3.4.1}. \\
  $(ii)$- Since $\Omega$ of class $C^{0,1}$. Then, by the same way in  \cite[Theorem 2.1]{3}, we can prove that $\Omega$ is a $W^{K,p(x,y)}$-extension domain. So, for any $u \in W^{K,p(x,y)}(\Omega)$ we define the extension function $\widetilde{u}$ by 
  $$ \widetilde{u}(x)= \left\{
   \begin{array}{clc}
   \ u(x)&if&x\in\Omega,\\\\
    ~~~~0 &if& x\in\mathbb{R}^{N}\setminus\Omega.
   \end{array}
   \right.$$
  The rest of proof is similar to Lemma \ref{3.4.1}. \hspace{8cm}$\Box$\\
  
  In the following lemma we prove a convergence property for a bounded sequence in $W^{K,p(x,y)}_{0}(\Omega)$. 
  
  \begin{lemma}
  Under the same assumptions of Theorem $\ref{3.1.1}$. And let $\{u_j\}$ be a bounded sequence in $W^{K,p(x,y)}_0(\Omega)$. Then there exists $u\in L^{r(x)}(\mathbb{R}^{N})$, with $u=0$ a.e in $\mathbb{R}^{N}\setminus \Omega$, such that up to a subsequence $$ u_n\longrightarrow u \hspace{0.3cm}strongly~~ in\hspace{0.3cm} L^{r(x)}(\Omega), \hspace{0.3cm}as\hspace{0.3cm} n\longrightarrow+\infty.$$
  \end{lemma}
\textbf{\textit{Proof}}. Since  $u_j\in W^{K,p(x,y)}_0(\Omega)$, then Lemma \ref{3.3}-$(ii)$ implies that $u_j\in W^{s,p(x,y)}(\mathbb{R}^{N})$, hence 
  $u_j\in W^{s,p(x,y)}(\Omega)$. Moreover, by Lemma \ref{3.3}-$(ii)$, Theorem \ref{3.1.1}-($ii$) and the definition of $W^{K,p(x,y)}_0(\Omega)$, we have
  $$\|u_j\|_{W^{s,p(x,y)}(\Omega)}\leqslant\|u_j\|_{W^{s,p(x,y)}(\mathbb{R}^N)}\leqslant C_0[u_j]_{K,p(x,y)}.$$
  Using this fact and since $\{u_j\}$ is bounded in $W^{K,p(x,y)}_0(\Omega)$, we get that $\{u_j\}$ is bounded in $ W^{s,p(x,y)}(\Omega)$. By Theorem \ref{1.1}, there exists $u\in L^{r(x)}(\Omega)$, such that up to a subsequence  $ u_n\longrightarrow u$ strongly in $L^{r(x)}(\Omega)$.
  Since $u_j=0$ a.e. in $\mathbb{R}^{N}\setminus \Omega$, we can define  $u=0$ a.e. in $\mathbb{R}^{N}\setminus \Omega$. \hspace{14.3cm}$\Box$\\ 
  
  As in the classic case with $s$ being an integer, any function in the
  fractional Sobolev space $W^{K,p(x,y)}(\Omega)$ can be approximated by a sequence of smooth functions with compact support.
  \begin{lemma}
  Let $(\ref{1})$, $(\ref{2})$ and $(\ref{3})$ be satisfied. Then the space $C^{\infty}_{0}(\mathbb{R}^{N})$ of smooth functions with
  compact support is  dense in $ W^{K,p(x,y)}(\Omega)$. 
  \end{lemma}
 \textbf{\textit{Proof}}. The proof is similar to the model case $K(x,y)=|x-y|^{-(N+sp(x,y))}$, in \cite[Lemma 2.3]{2}.\hspace{14cm}$\Box$
  
  \begin{remark}
  It is worth mentioning that our functional setting above is inspired by the pioneering works of M. Xiang et al. in \cite{20} when  \small$1<p(x,y)=p=constant<+\infty$, Servadi and Valdinoci in \cite{18,19} in which the corresponding functional framework was discussed as $p=2$.
  \end{remark}
  \section{Properties of the nonlocal fractional operator $\mathcal{L}^{p(x,.)}_K$}\label{400}
  
  In this section we give some basic properties of the nonlocal integro-differential operator of elliptic type  $\mathcal{L}^{p(x,.)}_K$.\\
  
  Let (\ref{1}) and (\ref{2}) be satisfied and    $K:\mathbb{R}^N \times\mathbb{R}^N\longrightarrow (0,+\infty)$ is a measurable function satisfy (\ref{4})-(\ref{6}). Then 
  $$\hspace{-3cm}\mathcal{L}^{p(x,.)}_K:W^{K,p(x,y)}_{0}(\Omega)\longrightarrow \left( W^{K,p(x,y)}_{0}(\Omega)\right)^{\ast}$$
  $$\hspace{2.4cm}u\longrightarrow\mathcal{L}^{p(x,.)}_K(u) : W^{K,p(x,y)}_{0}(\Omega)\longrightarrow \mathbb{R}$$
  $$\hspace{9.6cm}\varphi\longrightarrow\mathcal<\mathcal{L}^{p(x,.)}_K(u),\varphi>$$
  such that 
  $$<\mathcal{L}^{p(x,.)}_K(u),\varphi>=\int_{\mathbb{R}^{2N}}|u(x)-u(y)|^{p(x,y)-2}(u(x)-u(y))(\varphi(x)-\varphi(y))K(x,y)~dxdy,$$
  where $\left( W^{K,p(x,y)}_{0}(\Omega)\right) ^{\ast}$ is the dual space of $W^{K,p(x,y)}_{0}(\Omega)$.\\
  
  In the following Lemma, we show some fundamental properties of the operator $\mathcal{L}^{p(x,.)}_K$.
  
  \begin{lemma}\label{4.1}
  Suppose that $(\ref{1})$ and $(\ref{2})$ be satisfied and let $K:\mathbb{R}^N \times\mathbb{R}^N\longrightarrow (0,+\infty)$ be a measurable function satisfying $(\ref{4})$-$(\ref{6})$. Then, The following assertions hold: 
   \begin{enumerate}[label=(\roman*)]
   \item $\mathcal{L}^{p(x,.)}_K$ is well defined and bounded,
   \item $\mathcal{L}^{p(x,.)}_K$ is a strictly monotone operator,
   \item  $\mathcal{L}^{p(x,.)}_K$ is a mapping of type $(S_{+})$, that is, if $u_{k}\rightharpoonup u$ in $W^{K,p(x,y)}_{0}$ and $\underset{k\longrightarrow +\infty}{\limsup}<\mathcal{L}^{p(x,.)}_K(u_{k})-\mathcal{L}^{p(x,.)}_K(u), u_{k}-u>\leqslant0$, then  $u_{k}\longrightarrow u$ in $W^{K,p(x,y)}_{0}$,
   \item $\mathcal{L}^{p(x,.)}_K:W^{K,p(x,y)}_{0}(\Omega)\longrightarrow \bigg(W^{K,p(x,y)}_{0}(\Omega)\bigg)^{\ast}$ is a homeomorphism,
   \item $\mathcal{L}^{p(x,.)}_K$ is coercive.
   \end{enumerate}
  \end{lemma} 
  \textbf{\textit{Proof}}.$(i)$- Let $u,\varphi \in W^{K,p(x,y)}_{0} (\Omega)$. Then,  \begingroup\makeatletter\def\f@size{9.6}\check@mathfonts$$\begin{array}{clll}
  \big|<\mathcal{L}^{p(x,.)}_K(u),\varphi>\big|&\leqslant&\bigg|\int_{\mathbb{R}^{2N}}|u(x)-u(y)|^{p(x,y)-2}(u(x)-u(y))(\varphi(x)-\varphi(y))K(x,y)~dxdy\bigg|\\\\&\leqslant& \int_{\mathbb{R}^{2N}}\bigg(|u(x)-u(y)|^{p(x,y)-1}K(x,y)^{\frac{1}{\hat{p}(x,y)}}\bigg)\bigg(|\varphi(x)-\varphi(y)|K(x,y)^{\frac{1}{p(x,y)}}\bigg)~dxdy,
  \end{array}$$\endgroup
  where $\hat{p}:\mathbb{R}^{N}\times\mathbb{R}^{N}\longrightarrow(1,+\infty)$ is the conjugate exponent of $p$, that is,  \begingroup\makeatletter\def\f@size{9.3}\check@mathfonts$\frac{1}{\hat{p}(x,y)}+ \frac{1}{p(x,y)}=1.$\endgroup
  \\
 If we set $$\varPsi(x,y)=|u(x)-u(y)|^{p(x,y)-1}K(x,y)^{\frac{1}{\hat{p}(x,y)}}\in L^{\hat{p}(x,y)}(Q,dxdy),$$
  $$\varPhi(x,y)=|\varphi(x)-\varphi(y)|K(x,y)^{\frac{1}{p(x,y)}}\in L^{p(x,y)}(Q,dxdy).$$
 So, by H\"{o}lder inequality, we obtain
   $$\begin{array}{clll}
    \big|<\mathcal{L}^{p(x,.)}_K(u),\varphi>\big|&\leqslant&2\|\varPsi\|_{L^{\hat{p}(x,y)}(Q,dxdy)} \|\varPhi\|_{L^{p(x,y)}(Q,dxdy)}
    \\\\&\leqslant& C \|\varPhi\|_{L^{p(x,y)}(Q,dxdy)}.
    \end{array}$$
    It follows that 
    $$\|\mathcal{L}^{p(x,.)}_K(u)\|_{\big(W^{K,p(x,y)}_{0}(\Omega)\big)^{\ast}}\leqslant C <+\infty.$$
    For the proof of the properties $(ii), (iii)$ and $(iv)$, we follow the same argument in Lemma 4.2-($(i), (ii)$ and $(iii)$) in \cite{2}.\\
    $(v)$- Let  $u\in W^{K,p(x,y)}_{0}(\Omega)$. Then, we have 
    $$<\mathcal{L}^{p(x,.)}_K(u),u>= \int_{\mathbb{R}^{2N}}|u(x)-u(y)|^{p(x,y)}K(x,y)~dxdy=\rho^{o}_{K,p(.,.)}(u).$$
    If $[u]_{K,p(x,y)}>1$. From Lemma  \ref{2.3}-$(i)$, we get 
    \begin{equation}\label{17}
    <\mathcal{L}^{p(x,.)}_K(u),u>=\rho^{o}_{K,p(.,.)}(u)\geqslant [u]_{K,p(x,y)}^{p^{-}}=\|u\|_{W^{K,p(x,y)}_{0}}^{p^{-}}.
    \end{equation}
     If $[u]_{K,p(x,y)}<1$. By Lemma \ref{2.3}-$(ii)$, we get 
        \begin{equation}\label{18}
        <\mathcal{L}^{p(x,.)}_K(u),u>=\rho^{o}_{K,p(.,.)}(u)\geqslant [u]_{K,p(x,y)}^{p^{+}}=\|u\|_{W^{K,p(x,y)}_{0}}^{p^{+}}.
        \end{equation}
        Combining (\ref{17}) and (\ref{18}), we obtain
        $$\lim_{\|u\|_{W^{K,p(x,y)}_{0}}\longrightarrow +\infty}\frac{ <\mathcal{L}^{p(x,.)}_K(u),u>}{\|u\|_{W^{K,p(x,y)}_{0}}}=+\infty.$$
        This concludes the proof.\hspace{10.5cm}$\Box$
        
        \section{Application to nonlocal fractional problems with variable exponent}\label{500}
        In this section, we work under the hypotheses of Theorem \ref{3.1.1}. we aim to study two problems driven by the nonlocal operator $\mathcal{L}^{p(x,.)}$ and its particular case $(-\Delta_{p(x,.)})^{s}$
        \subsection{Application to Kirchhoff type problems} 
        At first, we discus the existence of a nontrivial solution for a nonlocal $p(x,.)$-Kirchhoff type problem of the following form 
        \begingroup\makeatletter\def\f@size{9.3}\check@mathfonts$$\label{01}\left( \mathcal{P}_{M}^{K}\right) ~~ \left\{
               \begin{array}{clccl}
                M \left( \ \int_{Q}\frac{|u(x)-u(y)|^{p(x,y)}}{p(x,y)}K(x,y)~dxdy\right) \mathcal{L}^{p(x,.)}_Ku(x)+ |u|^{\bar{p}(x)-2}u &=&f(x,u)&\text{in}&~~~\Omega,\\ \\
             \hspace{7cm}  u&=&\hspace{-0.3cm}0  &\text{in}&~~ ~\mathbb{R}^{N}\setminus\Omega, 
               \end{array}
               \right.
               $$   
        \endgroup
       where $\Omega\subset \mathbb{R}^{N}$, $N\geqslant3$,  is a Lipschitz bounded open domain, $M:\mathbb{R}^{+}\longrightarrow \mathbb{R}$ is a continuous function which satisfies the following polynomial growth condition\\
       
       \label{M} $\left(M_1\right) $ : $(1-\mu)t^{\alpha(x)-1}\leqslant M(t)\leqslant(1+\mu)t^{\alpha(x)-1}$, for all $t>0$ and $\mu\in[0,1),$\\with   $\alpha:\overline{\Omega}\longrightarrow(1,+\infty)$ is a bounded function such that $1<\alpha^-\leqslant \alpha(x)\leqslant \alpha^+<\infty.$
       $f:\Omega\times\mathbb{R}\longrightarrow\mathbb{R}$ is a Carath\'{e}odory function satisfies the following growth condition
          $$\label{f0}(f_{0}):\hspace{1cm}|f(x,t)|\leqslant c_1(1+|t|^{\beta(x)-1})\quad\text{ for all } (x,t)\in \Omega\times\mathbb{R}, $$
          where $\beta\in C_+(\overline{\Omega})$ such that $\beta(x)<p^{*}_s(x)$ for all $x\in \overline{\Omega}$, and $\frac{\beta^-}{\alpha^+}>p^+$.
          
        \label{f1}$(f_{1})$ : $\lim _{t\rightarrow0}\frac{f(x,t)}{|t|^{p^{+}-1}}=0$ \quad uniformly for $x\in \Omega.$ \\
          \label{AR}$(AR)$ :  There exist $A>0$ and $\theta> \left( \frac{1+\mu}{1-\mu}\right)  \frac{\alpha^+(p^+)^{\alpha^+}}{(p^-)^{\alpha^- -1}} $ such that $$0<\theta F(x,t)=\theta\int_{0}^{t}f(x,\tau)d\tau\leqslant f(x,t)t\quad \text{ for all } |t|>A \text{ and } \text{ a.e. } x\in \Omega.$$
            Actually, Ambrosetti-Rabinowitz condition \hyperref[AR]{$(AR)$} is quite natural and important not only to ensure that the Euler-Lagrange functional has a mountain pass geometry, but also to guarantee that the boundedness of Palais-Smale (PS) sequences.\\
            
            One typical feature of problem \hyperref[01]{$(\mathcal{P}_{M}^{K})$} is the nonlocality, in the sense that the value of $(-\Delta_{p(x,.)})^{s}u(x)$ at any point $x\in \Omega$ depends not only on the values of $u$ on $\Omega$, but actually on the entire space $\mathbb{R}^{N}$.   Moreover, the presence of the function $M$, which
            implies that the first equation in \hyperref[01]{$(\mathcal{P}_{M}^{K})$} is no longer a pointwise equation, it  is no longer a pointwise identity, therefore it is often called nonlocal problem.  Therefore, the Dirichlet datum is given in $\mathbb{R}^{N}\setminus\Omega$ (which is different from the classical case of the $p(x)$-Laplacian) and not simply on $\partial\Omega$. This causes
            some mathematical difficulties which make the study of such a problem particularly
            interesting. Motivated by the results in \cite{Af1,Af2,SM,SH}, we will prove that  problem \hyperref[01]{$(\mathcal{P}_{M}^{K})$} has at least one nontrivial weak solution, by means of mountain pass theorem of of Ambrosetti and Rabinowitz \cite{AR1}.\\ Throughout this part, for simplicity,  we use $c_i$, to denote the general nonnegative or positive constant
             (the exact value may change from line to line), we set also $X_0=W_0^{K,p(x,y)}(\Omega)$.
              \begin{definition}
              We say that $u\in X_0$ is a weak solution of problem \hyperref[01]{$(\mathcal{P}_{M}^{K})$} if
              $$  
             M\left( \ \sigma_{p(x,y)}(u) \right) \int_{Q}|u(x)-u(y)|^{p(x,y)-2}\big(u(x)-u(y)\big)\big (\varphi(x)-\varphi(y)\big)K(x,y)~dxdy~ 
              $$
                \begin{equation}\label{51}
             +\int_{\Omega} |u|^{\bar{p}(x)-2}u \varphi dx -\int_{\Omega}f(x,u)\varphi(x)dx=0,  
                \end{equation}
               for all $\varphi \in X_{0}$, where $$ \sigma_{p(x,y)}(u)=\int_{Q}\frac{|u(x)-u(y)|^{p(x,y)}}{p(x,y)}K(x,y)~dxdy~.$$
              \end{definition}
              
              Let us consider the Euler-Lagrange functional $J: X_{0}\longrightarrow \mathbb{R}$ which associated to  \hyperref[01]{$(\mathcal{P}_{M}^{K})$}, and defined by 
              $$\begin{array}{clc}
              J(u)&=&\widehat{M} \left( \ \int_{Q}\frac{|u(x)-u(y)|^{p(x,y)}}{p(x,y)}K(x,y)~dxdy~\right)+ \int_{\Omega}\frac{1}{\bar{p}(x)}|u|^{\bar{p}(x)}dx-\int_{\Omega}F(x,u)dx\\&=&\hspace{-4.5cm}\widehat{M} \left( \sigma_{p(x,y)}(u)\right)+ \int_{\Omega}\frac{1}{\bar{p}(x)}|u|^{\bar{p}(x)}dx-\int_{\Omega}F(x,u)dx,
              \end{array}$$
              where $\widehat{M}(t)=\int_{0}^{t}M(\tau)d\tau $.\\
              Standard arguments (see, for instance \cite[Lemma 3.1]{SH}) and the continuity of $M$ imply that $J$ is well defined and $J\in C^{1}(X_0,\mathbb{R})$. Moreover, for all $u,\varphi \in X_{0}$, its G\^{a}teaux derivative is given by 
             \begingroup\makeatletter\def\f@size{10}\check@mathfonts $$<J'(u),\varphi>= 
              M\left( \ \sigma_{p(x,y)}(u)\right) \int_{Q}|u(x)-u(y)|^{p(x,y)-2}\big(u(x)-u(y)\big)\big (\varphi(x)-\varphi(y)\big)K(x,y)~dxdy~ 
               $$
                 \begin{equation*}
              +\int_{\Omega} |u|^{\bar{p}(x)-2}u \varphi dx -\int_{\Omega}f(x,u)\varphi dx.   \end{equation*}\endgroup
               Thus, the weak solutions of \hyperref[01]{$(\mathcal{P}_{M}^{K})$} coincide with the critical points of $J$.\\
               
             Now, we are in a position to state our existence result as follows
             \begin{theorem}\label{t3.1}
             Let $\Omega$ be a Lipschitz bounded  domain in $\mathbb{R}^{N}$ and let $s\in (0,1)$, let $p:\overline{Q}\longrightarrow(1,+\infty)$  be a continuous function satisfies $(\ref{1})$ and  $(\ref{2})$ with $s p^{+}<N$. Assume that the assumptions \hyperref[M]{$(M_1)$}, \hyperref[f0]{$(f_0)$}, \hyperref[f1]{$(f_1)$} and \hyperref[AR]{$(AR)$} hold. Then, problem \hyperref[01]{$(\mathcal{P}_{M}^{K})$} has at least one nontrivial weak solution.
             \end{theorem}
              The proof of Theorem \ref{t3.1} based on mountain pass theorem of  Ambrosetti and Rabinowitz, and it follows from the following Lemmas. \begin{lemma}\label{l4}
              Suppose that the assumptions   \hyperref[M]{$(M_1)$}, \hyperref[f0]{$(f_0)$},  and \hyperref[AR]{$(AR)$} hold. Then, $J$ satisfies the (PS) condition.
              \end{lemma} 
              \textit{\textbf{Proof.}} Let us assume that there exists a sequence $\{u_n\}\subset X_0$ sch that \begin{equation}\label{c7}
              \left\{
                 \begin{array}{clc}
                 |J(u_n)|\leqslant c_2,\\\\
                 J'(u_n)\underset{n\rightarrow +\infty}{\longrightarrow}0.
                 \end{array}
                 \right.
              \end{equation} 
              Using \hyperref[M]{$(M_1)$}, \hyperref[AR]{$(AR)$}, Proposition \ref{2.2.1}, Lemma \ref{2.3} and  Remark \ref{r3.4}-$(i)$, for $n$ large enough, we get\\\\
              $c_2+\|u_n\|_{W_0}$
              $$\begin{array}{cll}
              &\geqslant& \hspace{-1cm}J(u_n)
              -\frac{1}{\theta}<J'(u_n),u_n>\\\\&=&\hspace{-1cm}\widehat{M}\left( \sigma_{p(x,y)}(u_n)\right)+\int_{\Omega}\frac{1}{\bar{p}(x)}|u|^{\bar{p}(x)}dx-\int_{\Omega}F(x,u_n)dx\\&\hspace{0.4cm}-& \hspace{-0.9cm}\frac{1}{\theta}M\left( \sigma_{p(x,y)}(u_n)\right)\int_{Q}|u_n(x)-u_n(y)|^{p(x,y)}K(x,y)~dxdy -\frac{1}{\theta}\int_{\Omega} |u_n|^{\bar{p}(x)} dx +\frac{1}{\theta}\int_{\Omega}f(x,u_n)u_ndx\\\\&\geqslant&\hspace{-1cm}\frac{1-\mu}{\alpha^+}\left( \sigma_{p(x,y)}(u_n)\right)^{\alpha(x)} -\frac{1+\mu}{\theta}\left( \sigma_{p(x,y)}(u_n)\right)^{\alpha(x)-1}\rho^0_{K,p(.,.)}(u_n)+\left( \frac{1}{p^+}-\frac{1}{\theta}\right) \int_{\Omega}|u_n|^{\bar{p}(x)}dx\\&\hspace{1cm}+&  \int_{\Omega} \left[\frac{1}{\theta}f(x,u_n)u_n-F(x,u_n)  \right] dx  
               \end{array}$$
             $$  \begin{array}{clc}
            &\geqslant&\frac{1-\mu}{\alpha^+(p^{+})^{\alpha^+}}\left(\rho_{K,p(.,.)}^0(u_n)\right)^{\alpha(x)} -\frac{1+\mu}{ \theta(p^{-})^{\alpha^--1}}\left( \rho^0_{K,p(.,.)}(u_n)\right)^{\alpha(x)}+\left( \frac{1}{p^+}-\frac{1}{\theta}\right)\|u_n\|_{L^{\bar{p}(x)}(\Omega)}^{p^{-}}\\\\&\geqslant&\hspace{-4.6cm} \left(\frac{1-\mu}{\alpha^+(p^{+})^{\alpha^+}}-\frac{1+\mu}{ \theta(p^{-})^{\alpha^--1}}\right) \|u_n\|_{X_0}^{\alpha^- p^{-}}+\left( \frac{1}{p^+}-\frac{1}{\theta}\right)c_3\|u_n\|_{X_0}^{p^{-}}.
               \end{array}$$
              Since $1<p^-<\alpha^- p^-$ and $\theta> \left( \frac{1+\mu}{1-\mu}\right)  \frac{\alpha^+(p^+)^{\alpha^+}}{(p^-)^{\alpha^- -1}}>p^+ $, we obtain that $\{u_n\}$ is bounded in $X_{0}$. This information, combined with the fact that $X_{0}$ is reflexive, implies that there exists a subsequence, still denote by $\{u_n\}$, and $u\in X_0 $ such that $\{u_n\}$ converges weakly to  $u$ in $X_0$. Next, as $\beta(x)<p^{*}_s(x)$ for all $x\in \overline{\Omega}$, then by Remark \ref{r3.4}-(1), $X_0$ is compactly embedded in $L^{\beta(x)}(\Omega)$,  it follows that \begin{equation}\label{5.3}
               u_n \longrightarrow u \text{ (strongly) in } L^{\beta(x)}(\Omega) \quad\text{and }  \quad u_n(x) \longrightarrow u(x) \text{ a.e. } x\in \Omega.
               \end{equation}
              Using (\ref{c7}), we have 
              $$<J'(u_n), u_n -u>\longrightarrow0,\quad\text{ as } n\rightarrow +\infty,$$
              that is,
              $$
               M\left( \ \sigma_{p(x,y)}(u_n)\right) <\mathcal{L}^{p(x,.)}_K(u_n),u_n-u> 
                +\int_{\Omega} |u_n|^{\bar{p}(x)-2}u_n(u_n-u) dx -\int_{\Omega}f(x,u_n)(u_n-u) dx \underset{n\rightarrow+\infty}{\longrightarrow}0   .$$
              From \hyperref[f0]{$(f_0)$} and Lemma \ref{2.1}, it follows that 
              $$\left| \int_{\Omega}f(x,u_n)(u_n-u)dx\right| \leqslant c \left\||u_n|^{\beta(x)-1} \right\|_{L^{\hat{\beta}(x)}(\Omega)}\|u_n-u\|_{L^{\beta(x)}(\Omega)}+ c_1\int_{\Omega}|u_n-u|dx,$$
              where $\frac{1}{\beta(x)}+\frac{1}{\hat{\beta}(x)}=1$. So, by (\ref{5.3}), we have 
               \begin{equation}\label{5.4}
               \int_{\Omega}f(x,u_n)(u_n-u)dx\longrightarrow 0 \quad\text{ as } n\rightarrow +\infty.
               \end{equation}
               Using again H\"{o}lder inequality,
we obtain              
        $$ \int_{\Omega} |u|^{\bar{p}(x)-2}u (u_n-u) dx\leqslant 2 \|u_n\|_{L^{\bar{p}(x)}(\Omega)} \|u_n-u\|_{L^{\bar{p}(x)}(\Omega)}.$$
       Since $u_n$ converges weakly to  $u$ in $X_0$, then  the compact embedding of $X_0$ into $L^{\bar{p}(x)}(\Omega)$ helps us to get 
        \begin{equation}\label{5.5}
       \lim_{n\rightarrow \infty}\int_{\Omega} |u_n|^{\bar{p}(x)-2}u (u_n-u) dx=0.
        \end{equation}
        Hence, by (\ref{5.4}) and (\ref{5.5}), we get
        \begin{equation}\label{5.6}
         M\left( \ \sigma_{p(x,y)}(u_n)\right) <\mathcal{L}^{p(x,y)}_K(u_n),u_n-u> 
          \underset{n\rightarrow+\infty}{\longrightarrow}0 . 
        \end{equation}
        Now, since $\{u_n\}$ is bounded in $X_0$, we may assume that $$\sigma_{p(x,y)}(u_n)\underset{n\rightarrow+\infty}{\longrightarrow}t_1\geqslant0.$$
        If $t_1=0$, then $\{u_n\}$ converge strongly to $u=0$ in $X_0$ and the proof is finished.\\
        If $t_1>0$, since the function $M$ is continuous, we have 
        
        $$M(\sigma_{p(x,y)}(u_n))\underset{n\rightarrow+\infty}{\longrightarrow}M(t_1)\geqslant0.$$
        Hence, by  \hyperref[M]{$(M_1)$}, for $n$ large enough, we have that \begin{equation}\label{5.7}
        0<c_4\leqslant M(\sigma_{p(x,y)}(u_n))\leqslant c_5.
        \end{equation}
        Combining (\ref{5.6}) and (\ref{5.7}), we deduce that 
        \begin{equation}\label{5.8}<\mathcal{L}^{p(x,y)}_K(u_n),u_n-u> 
          \underset{n\rightarrow+\infty}{\longrightarrow}0.\end{equation}
          On the other hand, since $\{u_n\}$ converge weakly to $u$ in $X_0$, we have that $$<J'(u),u_n-u> 
                    \underset{n\rightarrow+\infty}{\longrightarrow}0,$$
            that is,
               $$
                             M\left( \ \sigma_{p(x,y)}(u)\right) <\mathcal{L}^{p(x,.)}_K(u),u_n-u> 
                              +\int_{\Omega} |u|^{\bar{p}(x)-2}u (u_n-u) dx -\int_{\Omega}f(x,u)(u_n-u) dx \underset{n\rightarrow+\infty}{\longrightarrow}0,   $$      
which implies by using the same argument as before that   \begin{equation}\label{5.9}<\mathcal{L}^{p(x,y)}_K(u),u_n-u> 
          \underset{n\rightarrow+\infty}{\longrightarrow}0.\end{equation} 
          Combining \textsc{(\ref{5.8})}, and    \textsc{(\ref{5.9})} we deduce that 
            $$\limsup_{n\rightarrow+\infty}<\mathcal{L}^{p(x,y)}_K(u_n)-\mathcal{L}^{p(x,y)}_K(u),u_n-u>\leqslant0.$$
        By Lemma \ref{4.1}-$(iii)$ $\mathcal{L}^{p(x,y)}_K$ is a mapping of type $(S_{+})$, thus  $$\left \{ \begin{array}{clc}
           \limsup_{n\rightarrow+\infty} <\mathcal{L}^{p(x,y)}_K(u_n)-\mathcal{L}^{p(x,y)}_K(u),u_n-u>\leqslant0,\\
           u_k\rightharpoonup u$ in $X_0,\\
           \mathcal{L} ~~\text{is a mapping of type}~~ (S_{+}).
           \end{array}
              \right.\hspace{0.5cm}\Rightarrow u_n\longrightarrow u \text{ (strongly) in } X_0 .$$ 
              Consequently, $J$ satisfies the (PS) condition.\hspace{7.5cm}$\Box$\\
                   
            The following lemma shows that the functional $J$
             satisfies the first geometrical condition of the mountain pass theorem;
             \begin{lemma}\label{l5.2}
            Suppose that the assumptions   \hyperref[M]{$(M_1)$}, \hyperref[f0]{$(f_0)$},  and \hyperref[f1]{$(f_1)$} hold. Then there exist two positive real numbers $R$ and $a$ such that $J(u)\geqslant a>0$ for all $u\in X_0$ with $\|u\|_{X_0}=R$. 
             \end{lemma}   
            \textit{\textbf{Proof.}} Let $u\in X_0$ with $\|u\|_{X_0}<1$, then by \hyperref[M]{$(M_1)$}, we have 
              \begin{equation}\label{5.10}
              J(u)\geqslant \frac{1-\mu}{\alpha^+(p^+)^{\alpha^+}}\|u\|_{X_0}^{\alpha^+p^+}+\int_{\Omega}\frac{1}{\bar{p}(x)}|u|^{\bar{p}(x)}dx-\int_{\Omega}F(x,u)dx.
              \end{equation}
              Since $\beta(x)<p^{*}_s(x)$  and $p^+<p^{*}_s(x)$ for all $x\in \overline{\Omega},$ then by  Remark \ref{r3.4}-(1), we have that $X_0$ is continuously embedded in $L^{\beta(x)}(\Omega)$, $L^{\bar{p}(x)}(\Omega)$ and $L^{p^+}(\Omega)$, that is , there exist three positive constants $c_6$, $c_7$ and $c_8$ such that  \begin{equation}\label{5.11}
                \|u\|_{L^{\beta(x)}(\Omega)}\leqslant c_6\|u\|_{X_0},\quad \|u\|_{L^{p^+}(\Omega)}\leqslant c_7\|u\|_{X_0}\quad\text{ and }\quad \|u\|_{L^{\bar{p}(x)}(\Omega)}\leqslant c_8\|u\|_{X_0}.
                \end{equation}
              Now, we assume that  $\|u\|_{X_0}<\min \left\lbrace 1,\frac{1}{c_6}, \frac{1}{c_8}\right\rbrace $, where $c_6$ and $c_8$ are the positive constant given in (\ref{5.11}), then we get 
                $$\|u\|_{L^{\bar{p}(x)}(\Omega)}<1\quad\text{and} \quad\|u\|_{L^{\beta(x)}(\Omega)}<1 \text{ for all } u\in X_0~~\text{ with } \|u\|_{X_0}=R\in (0,1).$$
                Combining  \hyperref[f0]{$(f_0)$} and  \hyperref[f1]{$(f_1)$}, we get  $$F(x,t)\leqslant \varepsilon |t|^{p^+}+c_{\varepsilon}|t|^{\beta(x)}\quad\text{ for all } (x,t)\in\Omega\times\mathbb{R}.$$
                Therefore, by \hyperref[M]{$(M_1)$} ,(\ref{5.10}), (\ref{5.11}) and Proposition \ref{2.2.1}-$(i)$, we obtain
               $$\begin{array}{clc}
                J(u)&\geqslant& \frac{1-\mu}{\alpha^+(p^+)^{\alpha^+}}\|u\|_{X_0}^{\alpha^+p^+}
             +\frac{1}{p^+}\int_{\Omega}|u|^{\bar{p}(x)}dx  -\varepsilon\int_{\Omega}|u|^{p^+}dx-c_{\varepsilon}\int_{\Omega}|u|^{\beta(x)}dx \\\\
              &\geqslant&\hspace{-1.3cm}\frac{1-\mu}{\alpha^+(p^+)^{\alpha^+}}\|u\|_{X_0}^{\alpha^+p^+}+\frac{1}{p^+}\|u\|^{p^+}_{L^{\bar{p}(x)}(\Omega)}
              -\varepsilon\|u\|^{p^+}_{L^{p^+}(\Omega)}- c_{\varepsilon}\|u\|^{\beta^-}_{L^{\beta(x)}(\Omega)}
              \\\\
              &\geqslant&\hspace{-2cm}\frac{1-\mu}{\alpha^+(p^+)^{\alpha^+}}\|u\|_{X_0}^{\alpha^+p^+}+\frac{c_8^{p_+}}{p^+}\|u\|^{p^+}_{X_0}
              -\varepsilon c_{7}^{p^+}\|u\|_{X_0}^{p^+}- c_{\varepsilon}c_6^{\beta^-}\|u\|_{X_0}^{\beta^-}
              \\\\
              &\geqslant&\hspace{-1.3cm}\left[ \frac{1-\mu}{\alpha^+(p^+)^{\alpha^+}}\|u\|_{X_0}^{\alpha^+p^+-p^+}
             +\frac{c_8^{p_+}}{p^+} -\varepsilon c_{7}^{p^+}- c_{\varepsilon}c_6^{\beta^-}\|u\|_{X_0}^{\beta^--p^+}\right]\|u\|_{X_0}^{p^+}. 
               \end{array}$$
               We introduce the function $g:[0,1]\longrightarrow\mathbb{R}$, defined by 
               $$g(t)=\frac{1-\mu}{\alpha^+(p^+)^{\alpha^+}}t^{\alpha^+p^+-p^+}
               - c_{\varepsilon}c_6^{\beta^-}t^{\beta^--p^+}.$$
               Since $\frac{\beta^-}{\alpha^+}>p^{+}$, it is clear that there exists $\bar{t}\in [0,1]$ such that $$\max_{t\in[0,1]}g(t)=g(\bar{t})>0.$$
               Hence, for a fixed $\varepsilon\in \left( 0,\frac{g(\bar{t})+\frac{c_8^{p_+}}{p^+} }{c_7^{p^+}}\right) $ small enough, there exist two positive real numbers $R$ and $a$ such that $J(u)\geqslant a>0$ for all $u\in W_0$ with $\|u\|_{X_0}=R\in(0,1)$. 
               \hspace*{1cm}$\Box$\\
               
                The following result shows that the functional $J$
                  satisfies the second geometrical condition of  mountain pass theorem;
                \begin{lemma}\label{l6}
                Assume that \hyperref[M]{$(M_1)$} and \hyperref[AR]{$(AR)$} hold. Then there exists $u_0\in X_0$   such that $\|u\|_{X_0}>R$, $J(u)<0.$
                \end{lemma}
                \textit{\textbf{Proof.}} From the assumption \hyperref[AR]{$(AR)$}, we have that $$F(x,tu)>t^{\theta}F(x,u) \quad\text{ for all } t\geqslant1 \text{ and a.e.} x\in \Omega. $$
                Hence, by \hyperref[M]{$(M_1)$}, for $v\in X_0$, $v\neq0$ and $t>1$, we have  
                $$\begin{array}{clc}
                J(tv)&=&\hspace{-2.6cm}\widehat{M} \left( \sigma_{p(x,y)}(tv)\right)+\int_{\Omega}\frac{1}{\bar{p}(x)}|tv|^{\bar{p}(x)}dx-\int_{\Omega}F(x,tv)dx\\&\leqslant&\frac{1+\mu}{ \alpha^-(p^{-})^{\alpha^-}}t^{\alpha^+p^+}(\rho^0_{K,p(.,.)}(v))^{\alpha^+}+\frac{1}{p^-}t^{p^+}\int_{\Omega}|v|^{\bar{p}(x)}dx-t^{\theta}\int_{\Omega}F(x,v)dx.
                \end{array}$$
               From $\theta> \left( \frac{1+\mu}{1-\mu}\right)  \frac{\alpha^+(p^+)^{\alpha^+}}{(p^-)^{\alpha^- -1}}$, we find that $\theta>\alpha^+p^+>p^+$. Therefore, $$J(tv)\longrightarrow-\infty\quad\text{as}\quad t\rightarrow+\infty.$$ \hspace*{14.5cm}$\Box$\\
               Now, we are ready to prove  Theorem \ref{t3.1}.\\
               
                \textit{\textbf{Proof} of Theorem \ref{t3.1}.} Combining Lemmas \ref{l4}-\ref{l6} and the fact that $J(0)=0,$  we have that $J$ satisfies the assumptions of mountain pass theorem (see \cite{AR1}). Therefore, $J$ has at least one nontrivial critical point, that is, problem \hyperref[01]{$(\mathcal{P}_{M}^{K})$} has at least one nontrivial weak solution. \hspace{13cm}$\Box$
                \begin{example}
            \textnormal{ As a particular case we can take\\
              $\bullet$ $ M(t)= a+bt^{\alpha(x)-1}, \quad a,b>0 $, with   $\alpha:\overline{\Omega}\longrightarrow(1,+\infty)$ is a bounded function such that $1<\alpha^-\leqslant \alpha(x)\leqslant \alpha^+<\infty.$\\
              $\bullet$ $K(x,y)=|x-y|^{-(N+sp(x,y))}.$\\
              $\bullet$ $f(x,t)=|t|^{\gamma(x)-2}t$, where $\gamma\in C_+(\overline{\Omega})$  such that $\gamma(x)<p_s^*(x)$ for all $x\in \overline{\Omega}$ and $\gamma^->\alpha^+ p^+$.\\
              In this case, problem \hyperref[01]{$(\mathcal{P}_{M}^K)$} becomes}
              \begingroup\makeatletter\def\f@size{8}\check@mathfonts$$\label{011}\left( \mathcal{P}_{a,b}^s\right) ~~ \left\{
                             \begin{array}{clccl}
                              \left( a+b \left( \int_{Q}\frac{|u(x)-u(y)|^{p(x,y)}}{p(x,y)|x-y|^{N+sp(x,y)}}~dxdy\right) ^{\alpha (x)-1}\right) (-\Delta_{p(x,.)})^{s}u(x )+ |u|^{\bar{p}(x)-2}u &=&|u|^{\gamma(x)-2}u&\text{in}&~~~\Omega,\\ \\
                           \hspace{7cm}  u&=&\hspace{-0.3cm}0  &\text{in}&~~ ~\mathbb{R}^{N}\setminus\Omega, 
                             \end{array}
                             \right.
                             $$   
                      \endgroup
                \end{example} 
               It is easy to see that  the function $M$ satisfies  \hyperref[M]{$(M_1)$} and the function $f$ verify the assumptions   \hyperref[f0]{$(f_0)$}, \hyperref[f1]{$(f_1)$} and the \hyperref[AR]{$(AR)$} conditions. Consequently, problem \hyperref[011]{$(\mathcal{P}_{a,b}^s)$} has at least one nontrivial weak solution.

        \subsection{Existence and uniqueness result for a nonlocal problems}
         Now, we investigate the existence of a unique weak solution for a variational problem driven by general integro-differential operators of nonlocal fractional type $\mathcal{L}^{p(x,.)}_{K}$.
         $$\label{001}(\mathcal{P}_{{K}})~~ \left\{
              \begin{array}{clclc}
              \mathcal{L}^{p(x,.)}_K(u(x)) &=&f &\text{in}& \Omega,\\ \\
             \hspace{1.7cm} u&=&0 ~~~~~~~~~~~~~~~~~~~~~~~~~~ &\text{in}& ~\mathbb{R}^{N}\setminus\Omega,
              \end{array}
              \right.
              $$
              where $\Omega$ is a bounded open set of $\mathbb{R}^{N}$ and $f\in  X_{0} ^{\ast} $. 
              \begin{definition}
              We say that $u\in X_{0}$ is a weak solution of problem \hyperref[001]{$(\mathcal{P}_{K})$}, if 
              \begin{equation}
              \int_{\mathbb{R}^{2N}}|u(x)-u(y)|^{p(x,y)-2}(u(x)-u(y))(\varphi(x)-\varphi(y))K(x,y)~dxdy= \int_{\Omega} f v dx,
              \end{equation}
              for any $\varphi \in X_{0}$.
              \end{definition}
              Applying the Minty-Browder theorem, we get the following existence result. 
              \begin{theorem}\label{5.1}
              Let  $\Omega$ be a bounded open set of  $\mathbb{R}^{N}$ and $p:\mathbb{R}^{N}\times\mathbb{R}^{N}\longrightarrow(1,+\infty)$ be a continuous variable exponent satisfies $(\ref{1})$ and  $(\ref{2})$ with $sp^+<N$.  Suppose that  $K:\mathbb{R}^N \times\mathbb{R}^N\longrightarrow (0,+\infty)$ is a measurable function satisfying $(\ref{4})$-$(\ref{6})$  and $f\in X_{0}^{\ast} $. Then the problem \hyperref[001]{$(\mathcal{P}_{K})$} has a unique weak solution $u\in  X_{0}$.
              \end{theorem}
              A typical example for $K$ is given by the singular kernel $K(x,y)=|x-y|^{-(N+sp(x,y))}$. In this case, problem \hyperref[001]{$(\mathcal{P}_{K})$} becomes:
                       $$\label{010}(\mathcal{P}_{{s}})~~ \left\{
                            \begin{array}{clclc}
                            (-\Delta_{p(x,.)})^{s}(u(x)) &=&f &\text{in}& \Omega,\\ \\
                           \hspace{1.7cm} u&=&0 ~~~~~~~~~~~~~~~~~~~~~~~~~~ &\text{in}& ~\mathbb{R}^{N}\setminus\Omega.
                            \end{array}
                            \right.
                            $$
                            
  As a particular case, we derive an existence result for  problem \hyperref[010]{$(\mathcal{P}_{K})$}, which is given by the following corollary.
   \begin{corollary}\label{c5.2}
  Let  $\Omega$ be a bounded open set of  $\mathbb{R}^{N}$ and $p:\mathbb{R}^{N}\times\mathbb{R}^{N}\longrightarrow(1,+\infty)$ be a continuous variable exponent satisfy $(\ref{1})$ and  $(\ref{2})$ with $sp^+<N$. Let $s\in(0,1)$ and $f\in X_0^{\ast} $. Then the following equation 
   \begin{equation}\label{20}
                \int_{\mathbb{R}^{2N}}\frac{|u(x)-u(y)|^{p(x,y)-2}(u(x)-u(y))(\varphi(x)-\varphi(y))}{|x-y|^{N+sp(x,y)}}~dxdy= \int_{\Omega} f v dx,
                \end{equation}
                has a unique solution $u\in  X_0$.
  \end{corollary}                 
 \begin{remark}
 We observe that $(\ref{20})$ represents the weak formulation of problem \hyperref[010]{$(\mathcal{P}_{s})$}.
 \end{remark}                           
\textbf{\textit{Proof}} \textit{of Theorem} $\ref{5.1}$.    By Lemma \ref{4.1} the operator $\mathcal{L}^{p(x,.)}_{K}$  satisfies the conditions of Minty-Browder Theorem, that is,
\begin{itemize}
\item From Lemma \ref{4.1}-$(i)$, $\mathcal{L}^{p(x,.)}_{K}$ is bounded, from $X_0$ into $X_0^{\ast} $ .
\item From Lemma \ref{4.1}-$(ii)$, $\mathcal{L}^{p(x,.)}_{K}$ is a strictly monotone operator.
\item From Lemma \ref{4.1}-$(iv)$, $\mathcal{L}^{p(x,.)}_{K}$ is a homeomorphism. Hence, $\mathcal{L}^{p(x,.)}_{K}$ is continuous. 
\item From Lemma \ref{4.1}-$(v)$, $\mathcal{L}^{p(x,.)}_{K}$ is coercive.
\end{itemize}  
Consequently, in the light of Minty-Browder theorem \cite[Theorem V.15]{5}, then there exists a unique weak solution $u\in X_0$ of  problem \hyperref[001]{$(\mathcal{P}_{K})$}. \hspace{3.5cm}$\Box$\\
 
\textbf{\textit{Proof}} \textit{of Corollary} $\ref{c5.2}$. It is a consequence of Theorem \ref{5.1}, by choosing 
            $$K(x,y)=|x-y|^{-(N+sp(x,y))},$$
            and by recalling that $X_0\subset W^{s,p(x,y)}_{0}(\Omega)$.\hspace{8cm}$\Box$

\end{document}